%
%
%
\input amssym.def
\input amssym.tex
\magnification=\magstep1
\nopagenumbers \headline={\tenrm\hfil --\folio--\hfil}
\def\tomb{\phantom{.}\hfill\vrule height.4true cm width.3true cm \par\smallskip\noindent}
\def\ker{{\rm Ker\,}}  \def\chop{\hfill\break}  
\def\f #1,#2.{\mathsurround=0pt \hbox{${#1\over #2}$}\mathsurround=5pt}
\def\s #1.{_{\smash{\lower2pt\hbox{\mathsurround=0pt $\scriptstyle #1$}}\mathsurround=5pt}}
\def\r{{\hbox{\mathsurround=0pt$\rm I\! R$\mathsurround=5pt}}}

\def\mapdownl #1;{\vcenter{\hbox{$\scriptstyle#1$}}\Big\downarrow}
\def\mapdownr #1;{\Big\downarrow\rlap{$\vcenter{\hbox{$\scriptstyle#1$}}$}}
\tolerance=1600 \mathsurround=5pt  \def\aut{{\rm Aut}\,}
\def\maprightu #1;{\smash{\mathop{\longrightarrow}\limits^{#1}}}
\def\maprightd #1;{\smash{\mathop{\longrightarrow}\limits_{#1}}}

\def\brc #1,#2.{\left\langle #1\,|\,#2\right\rangle}
\def\rn#1{{\romannumeral#1}}  \def\cl #1.{{\cal #1}}
\def\convl #1,#2.{\mathrel{\mathop{\longrightarrow}\limits^{#1}_{#2}}}
\def\convr #1,#2.{\mathrel{\mathop{\longleftarrow}\limits^{#1}_{#2}}}
\def\set #1,#2.{\left\{\,#1\;\bigm|\;#2\,\right\}}
\def\theorem#1@#2@#3\par{\smallskip\parindent=.6true in \itemitem{\bf #1}
{\sl #2}\parindent=20pt\smallskip\itemitem{\it Proof:\/}#3\tomb}
\def\thrm#1"#2"#3\par{\smallskip\parindent=.6true in \itemitem{\bf #1}
{\sl #2}\parindent=20pt\smallskip\itemitem{\it Proof:\/}#3\tomb}
\def\teorem#1@#2@ {\smallskip\parindent=.6true in \itemitem{\bf Theorem #1}
{\sl #2}\hfill \parindent=20pt\smallskip\noindent} \def\ab{\allowbreak}
\def\fteorem#1@#2@ {\smallskip\parindent=.6true in \itemitem{\bf #1}{\sl #2}
\hfill\parindent=20pt\smallskip\noindent}

\def\shift #1;{\mathord{\phantom{#1}}}  
\def\ref #1.{\mathsurround=0pt${}^{#1}\phantom{|}$\mathsurround=5pt}
\def\cross #1.{\mathrel{\jackup 3,\mathop\times\limits_{#1}.}}
\def\jackup#1,#2.{\raise#1pt\hbox{\mathsurround=0pt $#2$\mathsurround=5pt}}

\def\bbrc #1,#2,#3.{\langle #1 |\,#2\,|#3\rangle}
   
\def\cst#1,{{C^*(#1-\un)}}
\def\un{\hbox{\mathsurround=0pt${\rm 1}\!\!{\rm 1}$\mathsurround=5pt}}

\def\alg#1.{{C^*(#1)}}

\def\ccr #1,#2.{{\overline{\Delta(#1,\,#2)}}}
\def\un{\Bbb I}
\def\r{\Bbb R}
\def\C{\Bbb C}
\def\b#1'{{\bf #1}}
\def\rest{\mathord{\restriction}}

\def\S{{\frak S}}

\def\supp{{\rm supp}\,}
\def\wt{\widetilde}
\def\wh{\widehat}
\def\wl{\wt{\cl L.}}

\pageno=1  \noindent
\centerline{\bf Host Algebras.}
\vglue .3in
\centerline{Hendrik Grundling}   
\centerline{ Department of Pure Mathematics, University of
New South Wales,}
\centerline{ Sydney, NSW 2052, Australia.}
\centerline{ email: hendrik@maths.unsw.edu.au}
\vglue .2in
\itemitem{{\bf Abstract}}{\sl  A host algebra generalises the concept of
a group algebra in the following way. Take a unital C*--algebra $\cl F.$
and a proper subset of its states $\S_0$ within which one wants to
keep the analysis
(e.g. the group algebra of a discrete group $G,$ and the set of its states
continuous w.r.t. some nondiscrete group topology of $G).$
Then a host algebra is a C*--algebra $\cl L.$ for which
we have 
embeddings  $\cl F.\subset\cl E.\supset\cl L.$
into a larger C*--algebra $\cl E.,$
such that states on $\cl L.$ extend uniquely to $\cl F.,$
and this extension defines a norm continuous affine bijection
between $\S_0$ and the whole state space of $\cl L..$
The main examples --though not the only ones--
are of course group and covariance algebras.
Here we study the general existence question for a host algebra of
a given pair $(\cl F.,\,\S_0),$ we show that given a host algebra one can do 
integral decompositions of states in $\S_0$ in terms of 
other states in $\S_0,$ and we show that if one does induction of
representations via host algebras, one stays within the class
of representations with the right continuity properties w.r.t.
$\S_0.$ Moreover, up to a central algebra, one can
always construct a host algebra (if $\S_0$ is a folium), 
but this central algebra can be an obstruction
to the existence of a host algebra.
These results should
be interesting to anyone who wants to construct a group
algebra for general topological groups, and quantum
physicists should also be interested due to various selection
criteria for physically acceptable states.}
\chop
\chop {\bf AMS classification: 46L05, 46L30, 81T05, 43A40}

\beginsection  Introduction.

In quantum physics one is frequently given a unital C*--algebra
$\cl F.$ for the observables and a distinguished proper subset of states 
$\S_0\subset\S(\cl F.)$ of its state space
together with the constraint that the physical system can only
realise these states and no others.
\itemitem{\bf Examples.}{\bf (1)} Constrained systems;-- here one has a
distinguished set of unitaries $\cl U.\subset\cl F._u$ and for
the set of physically realisable states $\S_0$ we have 
 the Dirac states $\set\omega\in\S(\cl F.),
{\omega(\cl U.)=1}.,$ cf. Grundling and Hurst [GH].
\itemitem{\bf (2)} The algebra of the canonical commutation
relations $\cl F.=\ccr S,B.$ over a symplectic space
$(S,\, B)$ (cf. Manuceau [Ma]), where $\S_0$ is required
to be the set of regular states.
\itemitem{\bf (3)} In finite dimensional quantum mechanics,
$\cl F.$ is taken as a factor of Type~I and the physically relevant states,
$\S_0,$ as its set of normal
states.
\itemitem{\bf (4)} In algebraic quantum field theory, let $\cl F.$
   be the inductive limit of a net of local algebras, and let 
   $\S_0$ be the set of locally normal states with respect
   to some distinguished state (cf. Haag [Ha]).

\noindent In such a situation one can object that the given system
$(\cl F.,\,\S_0)$ is not satisfactory because it leads naturally
to nonphysical objects, for instance the weak*--closure of
$\S_0$ can be the full state space $\S (\cl F.).$ Indeed, 
from the point of view that the states and observables
should be in some kind of duality (a Heisenberg--Schr\"odinger
picture equivalence), one can argue that if $\S_0$ is the
physical state space, then $\cl F.$ is not the correct algebra of
observables. From the mathematical point of view there are also problems,
e.g. we may not have a decomposition theory of states in $\S_0$
in terms of other states in $\S_0.$ Another problem, is that if
we have two pairs $(\cl F.^{(i)},\,\S_0^{(i)}),$ $i=1,\; 2$ and 
an imprimitivity bimodule for the algebras, then when we induce
a representation from one algebra to the other, we can easily move
out of the class of allowed states.
Our idea here is to replace $\cl F.$ by an algebra $\cl L.$
which has precisely $\S_0$ as its state space, in a sense to be 
made precise below. 

For some examples of $(\cl F.,\,\S_0),$ we do in fact have a more
convenient algebra which in some sense has precisely $\S_0$
as its state space, e.g. if we take $\cl F.=\cl B.(\cl H.),$ and let
$\S_0$ be its set of normal states, then the algebra of compact operators
$\cl L.=\cl K.(\cl H.)$ is  an algebra with state space $\S_0$
in the sense that $\S_0\restriction\cl K.(\cl H.)=\S(\cl K.(\cl H.))$
and the states on $\cl K.(\cl H.)$ extend uniquely.
(We use here the notation $\S(\cdot)$ for the state space of
its argument). Group algebras provide another set of examples 
-- see the next section.
So, inspired by these examples, we will study here the following
situation. Given a pair $(\cl F.,\,\S_0),$ find a pair of  C*--algebras
$\cl L.\subset\cl E.$ and an embedding $\xi:\cl F.\to\cl E.$
such that the states on $\cl L.$ extend uniquely to $\xi(\cl F.),$ 
and $\theta(\S(\cl L.))=\S_0$ where $\theta$ denotes the extension map
$\theta:\S(\cl L.)\to\S(\cl F.)$.
Naturally, there are existence questions to be answered, and we will
address these first.
In the next section we develop our basic theory, in Sect.~2 we
construct a host up to a central algebra, and in Sect.~3 we do 
a few applications.

\beginsection 1. Basic Concepts.

To prepare the ground, we first recall some background material.\chop
Recall that a hereditary subalgebra $\cl B.$ of a C*--algebra $\cl A.$
is a C*--subalgebra such that $0<A<B$ for $A\in \cl A.,$ $B\in\cl B.$
implies that $A\in\cl B.$. Equivalently (cf. Murphy [Mu]) $\cl B.$
is hereditary if $\cl BAB.\subseteq\cl B..$
Such algebras are plentiful, and are in bijection with the set of
closed left ideals of $\cl A..$ All closed two--sided ideals are hereditary. 
For us, the most important property is:
a C*--subalgebra $\cl B.$ is hereditary iff each state $\omega\in\S(\cl B.)$
has a unique extension to a state $\theta(\omega)\in\S(\cl A.)$
(cf. Kusuda [Ku]).
If we define a projection $P\in\cl A.''$ as the unit of 
$\cl B.''\subset\cl A.'',$ then the map $\theta:\cl B.^*\to
\cl A.^*$ defined by
$$\theta(\varphi)(A):=\varphi(PAP)=\lim_\alpha\varphi(E_\alpha AE_\alpha)
\eqno{(1)}$$
for all $A\in \cl A.$ and some approximate identity $\{E_{\alpha}\}$
of $\cl B.,$ is precisely the unique extension map on the states
$\S(\cl B.)\subset\cl B.^*.$ Moreover $\theta$ is an isometry, hence its range
is norm closed.
Given a representation of $\cl B.,$ we can always induce a representation
on $\cl A.$ from it (cf. Fell and Doran, XI.7.6~[FD]),
but usually this will be on a different space than the original representation.
In the case that $\cl B.$ is a two--sided ideal of  $\cl A.$
we have that $P\in\cl A.'\cap\cl A.''$, and so
$$\theta(\varphi)(A):=\varphi(PA)=\lim_\alpha\varphi(E_\alpha A)
\eqno{(2)}$$
and now, in addition, representations of $\cl B.$ also extend uniquely on
the same space to $\cl A..$ For a representation
$\pi:\cl B.\to\cl B.(\cl H.),$ the unique extension is
$$\wt\pi(A):=\pi(PA)=\mathop{\hbox{s-lim}}_\alpha\pi(E_\alpha A)
\eqno{(3)}$$
In this paper we will always use the notation $(\cl F.,\,\S_0),$
to denote a unital C*--algebra $\cl F.$ and a distinguished
subset of its state space. We choose $\cl F.$ to be unital,
since this ensures that its state space is w*--closed, hence
norm closed (cf. Pedersen [Pe] 3.2.1).
We are now ready for our basic definitions:
\item{\bf Def.} Given a pair $(\cl F.,\,\S_0)$ consisting of
a unital C*--algebra and a proper subset of its states,
we say a C*--algebra $\cl L.$ is a {\it host} for 
the pair if there is a unital C*--algebra $\cl E.\supset\cl L.$
(faithful embedding) and a unital *-homomorphism
${\xi:\cl F.\to\cl E.}$ such that:
\item{(\rn1)} $\cl L.$ is hereditary in $\cl E.,$
\item{(\rn2)} $\cl E.$ is generated by $\cl L.$ and $\xi(\cl F.),$
\item{(\rn3)} the map $\theta:\S(\cl L.)\to\S(\cl F.)$
defined by unique extension
$$\theta(\omega)(F)=\lim_\alpha\omega\big(E_\alpha\xi(F)E_\alpha)\;,\quad
F\in\cl F.\eqno{(4)}$$
for any approximate identity $\{E_{\alpha}\}$ of  $\cl L.,$ is injective 
and has range $\theta(\S(\cl L.))=\S_0.$
\item{$\bullet$} In the case that $\cl L.$ is in addition a two sided ideal
of $\cl E.,$ we call it an {\it ideal host}.  
\item{\bf Remarks} {\bf (1)} 
The requirement (\rn2) that  $\cl E.$ be generated by $\cl L.$ and $\xi(\cl F.),$
is not essential;-- if we start from some larger algebra $\cl E.$
satisfying the other requirements, we can always replace $\cl E.$
by $C^*(\cl L.\cup\xi(\cl F.))$ inside $\cl E.$
because $\cl L.$ is hereditary for any subalgebra of  $\cl E.$ containing it.
An important point in the definition, is that we allow $\cl L.$
to be outside $\cl F..$
\item{\bf (2)} If $(\cl F.,\,\S_0)$ has an ideal host $\cl L.\subset\cl E.,$
then there is a natural homomorphism of $\cl E.$ into the multiplier algebra
$M(\cl L.)\subset\cl L.'',$ and so we can equivalently define an ideal
host for $(\cl F.,\,\S_0)$ as a C*--algebra $\cl L.$ together with a homomorphism
$\xi:\cl F.\to M(\cl L.)$ such that $\theta:\S(\cl L.)\to\S(\cl F.)$
is injective, with image $\S_0.$
\item{\bf (3)} If $(\cl F.,\,\S_0)$ has only a host $\cl L.\subset\cl E.,$
then we may want to try a similar embedding than in the last remark.
The point is that $\cl E.$ acts as a set of quasi--multipliers of $\cl L.$
(i.e. $\cl L.E\cl L.\subset\cl L.$ for each $E\in\cl E.$)
and we know by C. Akemann and  G. Pedersen, 
Prop. 4.2 [AP],
 that there is a linear bijection between the quasi--multipliers 
and elements $A\in\cl L.''$ such that $\cl L.A\cl L.\subset\cl L..$
However, in general this bijection is not a homomorphism for
C*--algebras of quasi--multipliers, so we can not exploit this 
bijection as in the previous remark.
\item{\bf (4)} Ideal hosts are of course more useful than hosts
because then one has also unique extensions of representations
on the same space. As we shall see however,
 their existence place strong structural restrictions on $\S_0.$
Most of our analysis here will concern ideal hosts.

\noindent First we do a list of examples, and then some structural analysis.
\itemitem{\bf Examples}{\bf (1)} For the pair $(\cl F.,\,\S_0)=
\big(\cl B.(\cl H.),\;\S_N)$ where $\S_N$ denotes the set of
normal states, an ideal host algebra is $\cl L.=\cl K.(\cl H.)$
with the identity map $\xi:\cl B.(\cl H.)\to\cl B.(\cl H.)
=M(\cl K.(\cl H.))=\cl E.$.
\itemitem{\bf (2)} An important example for physics, is the following.
Let the pair $(\cl F.,\,\S_0)$ consists of the
CCR--algebra over a finite dimensional symplectic space,
and its set of regular states. To be more concrete, consider
 the CCR--algebra on $\r^2,$ i.e. the unique simple C*--algebra
$\cl F.$ generated by unitaries $\set\delta_{\bf x},{\bf x}\in\r^2.$
satisfying the Weyl relations:
$$\delta_{\bf x}\delta_{\bf y}=\rho({\bf x},\,{\bf y})\delta_{{\bf x}+{\bf y}}
\;,\qquad\hbox{where}\qquad
\rho({\bf x},\,{\bf y}):=\exp[i(x_1y_2-x_2y_1)]\,.$$
Define another C*--algebra $\cl L.$ as the C*--envelope
 of the twisted convolution
algebra, where the latter
 consists of $L^1(\r^2)$ equipped with the multiplication and involution:
$$f*g({\bf x})=\int_{\r^2}f({\bf y})\, g({\bf x}-{\bf y})\,\rho({\bf y},\,{\bf x})\,d{\bf y}\,,
\qquad\; f^*({\bf x})=\overline{f(-{\bf x})}.$$
This algebra $\cl L.$ is known to be isomorphic to $\cl K.(L^2(\r))$
(cf. I.E. Segal
 [Se]). Then $\cl F.\subset M(\cl L.)$ by the action 
$\delta_{\bf x}\cdot  f({\bf y})= \rho({\bf x},\,{\bf y})\, f({\bf y}-{\bf x}).$
 The unique extensions of states on $\cl L.$ to $\cl F.$ produce
precisely the set of regular states 
$$\S_0:=\set\omega\in\S(\cl F.),{\bf x}\to\omega(\delta_{\bf x})
\quad\hbox{is continuous}.\,.$$
The extension map $\theta:\S(\cl L.)\to\S_0$ is injective,
since for each $\phi\in\S_0$ we can reconstruct the $\omega\in\S(\cl L.)$
such that $\theta(\omega)=\phi$ via the formula
$$\omega(f):=\int_{\r^2}f({\bf x})\,\varphi(\delta_{{\bf x}})\,d{\bf x}$$
for all $f\in L^1(\r^2).$
(Put in other words, here $\cl F.$ is the twisted discrete group algebra
for $\r^2,$ and $\cl L.$ is the usual twisted group algebra
for $\r^2.)$
Since $\cl F.$ is simple, and $\S_0\not=\S(\cl F.),$ 
we see that $\cl F.\cap\cl L.=\{0\}.$ Of course $\cl L.$ is a far
better behaved C*--algebra than $\cl F.,$ it is even separable.
 See [Gr] for generalisations of this example.
\def\f{{\cl F.}}  \def\l{{\cl L.}} \def\b{{\cl B.}}

\itemitem{\bf (3)}
 Let $G$ be a nondiscrete topological group,
and denote $G$ equipped with the discrete topology by $G_d.$
Let the pair $(\cl F.,\,\S_0)$ be the discrete group algebra
$\cl F.=C^*(G_d)$ and its set of states continuous with respect
to the topology of $G$ (i.e. $g\to\omega(\delta_g)$ is continuous,
where $\delta_g$ denotes the Dirac point measure at $g$).
When $G$ is locally compact,
the usual group algebra $C^*(G)$ is an ideal host algebra
if we use the imbedding $C^*(G_d)\subset M\big(C^*(G)\big)$
obtained by convolution of measures.

\noindent Inspired by this example, we define
\item{\bf Def.} Let $G$ be a topological group, not necessarily locally compact.
Then a {\it group algebra} for it, is any ideal host for the pair
$(C^*(G_d),\,\S_0)$ where $\S_0$ denotes the states $\omega$ on $C^*(G_d)$
such that the map $g\to\omega(\delta_g)$ is continuous.

\noindent
In [Gr] we took this definition for a group algebra in order 
to construct a group algebra for
inductive limit groups. We can also adapt it for covariance algebras.
In this paper, however, we will not construct any such group algebras for
more general groups.

Next we analyze some structural consequences for the existence of a host for the
pair  $(\cl F.,\,\S_0).$ We will usually consider the homomorphism $\xi:
\cl F.\to\cl E.$ as
an embedding, to save on notation. 

\thrm Theorem 1.1." If a pair $(\cl F.,\,\S_0)$ has a host 
$\cl L.\subset\cl E.,$ then \chop
(\rn1) $\S_0$ is a norm--closed face in $\S(\cl F.).$\chop
(\rn2) If  $\cl L.$ is an ideal host, then the norm--closed face $\S_0$ is also
invariant, i.e. if $\omega\in\S_0$ then $\omega_B\in\S_0$
for all $B\in\cl F.$ with $\omega(B^*B)=1,$ and where 
$\omega_B(F):=\omega(B^*FB)$ for $F\in\cl F..$\chop
(\rn3) $\theta:\S(\cl L.)\to\S_0$ is an isomorphism, 
i.e. it is affine and a homeomorphism
w.r.t. the norm topology.
\chop (\rn4)
In the case that $\S_0$ is the face obtained by extending the states from
a hereditary subalgebra $\cl A.\subset\cl F.$
to $\cl F.,$ then ${\cl L.\cap\cl F.}=\cl A..$"
(\rn1) Norm closure: We first show that $\theta:\cl L.^*\to\cl F.^*$ is norm continuous.
Recall that since $\cl L.$ is hereditary in $\cl E.$ we have
$\theta(\varphi)(F):=\varphi(PFP)$ with $P\in\cl L.''\subset\cl E.''\supset\cl F.''.$
Now 
$$\eqalignno{\big\|\theta(\varphi)\big\|&=\sup\set\big|\theta(\varphi)(F)\big|,
{F\in\cl F.,\;\|F\|\leq 1}.\cr
&=\sup\set\big|\varphi(PFP)\big|,{F\in\cl F.,\;\|F\|\leq 1}.\cr}$$
Now we know from Kusuda [Ku] Theorem 2.2 that $P\cl E.''P=\cl L.''$
since $\cl L.$ is hereditary, and so
$$\set\big|\varphi(PFP)\big|,{F\in\cl F.,\;\|F\|\leq 1}.\subseteq
\set\big|\varphi(A)\big|,{A\in\cl L.'',\;\|A\|\leq 1}.  $$
and hence, since the supremum of the last set is just $\|\varphi\|,$
we find that $\big\|\theta(\varphi)\big\|\leq\|\varphi\|.$
Thus $\theta:\cl L.^*\to\cl F.^*$ is norm continuous.
By assumption $\theta$ is injective on $\S(\cl L.);$ we prove that
it is also injective on $\cl L.^*.$ If it were not, there would be
$\varphi_i\in\cl L.^*$ such that $\theta(\varphi_1-\varphi_2)=0.$
Then for $\psi:=\varphi_1-\varphi_2$ do a Jordan decomposition,
$\psi=\rho_+-\rho_-+i(\mu_+-\mu_-)$ and then $\theta(\psi)=0$
implies $\theta(\rho_+-\rho_-)=0=\theta(\mu_+-\mu_-),$ i.e.
$\theta(\rho_+)=\theta(\rho_-)$ and $\theta(\mu_+)=\theta(\mu_-).$
But $\theta$ is injective on states, so $\rho_+=\rho_-$ and $\mu_+=\mu_-,$
i.e. $\psi=0.$ Thus we know that $\theta$ is both norm continuous
and invertible on $\cl L.^*$ so by a corollary to the open mapping
theorem, its inverse must also be norm continuous, and hence by
Theorem 5.8, p216 of A. Taylor
[Ta] we conclude that $\theta(\cl L.^*)$ is norm closed.
Since $\S(\cl F.)\subset\cl F.^*$ is also norm closed (recall
that $\cl F.$ is unital), so is $\theta(\cl L.^*)\cap\S(\cl F.)$
and as this is just the image under $\theta$ of 
$\S(\cl L.),$ we conclude that $\S_0$ is norm closed.
Now (\rn3) also follows from the preceding.\chop
That $\S_0$ is a convex set, follows from the fact that $\theta$ is linear
and $\S(\cl L.)$ is convex, so we just need to prove the facial property.
Let $\omega=\lambda\varphi+(1-\lambda)\psi\in\S_0$ where $\lambda\in[0,1].$
We need to show that $\varphi,\;\psi\in\S_0.$ Since  $\omega\in\S_0$
there is a unique $\omega'\in\S(\cl L.)$ such that $\omega=\theta(\omega').$
By the hereditary property, $\omega'$ extends uniquely to a state 
$\wh{\omega'}$ on $\cl E.''$
(and this extension of course restricts to $\omega$ on $\cl F.$).
By definition of $P$ we have that $\wh{\omega'}(P)=1,$ and conversely, given any state 
$\gamma$ on $\cl E.''$ with $\gamma(P)=1$ we have that
$\gamma\restriction\cl L.\in\S(\cl L.),$ hence ${\theta(\gamma\restriction\cl L.)
}=\gamma\restriction\cl F.\in\S_0.$
Now $0=\wh{\omega'}(\un-P)=\lambda\wh\varphi(\un-P)+(1-\lambda)\wh\psi(\un-P)$
where $\wh\varphi,\;\wh\psi$ are extensions of $\varphi,\;\psi$
to $\cl E.''.$
Thus by positivity of all terms in the sum, we conclude that $\wh\varphi,\;\wh\psi$
vanish on $\un-P$ and hence $\varphi=\theta(\wh\varphi\restriction\cl L.),$
$\psi=\theta(\wh\psi\restriction\cl L.).$
Thus $\varphi,\;\psi\in\S_0,$ i.e. $\S_0$ is a face.\chop
For (\rn2), assume that $\cl L.$ is an ideal of $\cl E.,$ then we want to prove
invariance of the face. Let $\omega=\theta(\omega')\in \S_0,$
then clearly $\omega'_B(\cdot):=\omega'(B^*\,\cdot\,B)$ defines a state on 
$\cl L.$ using the fact that $\cl L.$ is an ideal (here we took $B\in\cl E.$
with $\omega'(B^*B)=1$). Thus $\S(\cl L.)$ is invariant, and by the definitions
$\omega_B=\theta(\omega'_B)\in\S_0$ because extension and conjugation of a state
commutes when the projection $P$ in Equation~(1) commutes with $\cl F.,$
and this it does since $\cl L.$ is an ideal. \chop
(\rn4) Since $\cl L.$ is hereditary in $\cl E.,$ we see that
 $\cl L.\cap\cl F.$ is hereditary in $\cl F..$ Moreover, the face 
$\S_0\subset\S(\cl F.)$ is precisely the states which extend from
the states on $\cl L.,$ which are uniquely determined by their values on
$\cl L.\cap\cl F.$ by the hereditary property. Since $\S_0$ are also
the states which extend from $\cl A.,$ by the bijection between 
such faces and hereditary subalgebras (cf. Pedersen 3.10.7 [Pe]
and Murphy 3.2.1 [Mu])
we conclude that $\cl L.\cap\cl F.=\cl A..$

\item{\bf Remarks.} {\bf (1)}
It is known that an invariant convex norm--closed set of states $\S_0$
is also a face (cf. footnote in [HKK]), so for an ideal host it suffices
to say that $\S_0$ is an invariant convex norm--closed set.
In other words, this says that the cone which  $\S_0$ generates, 
$\r_+\S_0,$ is a folium. We will extend the term ``folium'' to also mean
invariant norm--closed convex sets of state spaces.
\item{\bf (2)}
Since $\theta:\S(\cl L.)\to\S_0$ is an affine bijection,
it restricts to a bijection between the pure states on $\cl L.$
and the extreme points of $\S_0.$ This immediately limits the class
of faces and folia for which hosts exist, because there are many folia
without extreme points, e.g. the folium of normal states of
$L^\infty(\r)$ (measures absolutely continuous w.r.t. the
Lebesgue measure). We will sharpen this observation below.
Note also that as $\theta$ involves both an extension and a restriction,
it need not a priori take pure states to pure states.
Since $\S(\cl L.)$ is generated as the w*--closed convex hull 
of its pure states, we can use $\theta$ to transfer this
weak*--topology to $\S_0$ (but this is different from
the weak*--topology of $\cl F.),$ to conclude that w.r.t.
this topology $\S_0$ is a compact convex set, hence via
Choquet theory, there are integral decompositions 
of states w.r.t. measures on $\S_0.$
In the next section we will exploit these decompositions.
It would be nice to have some intrinsic definition
of the w*--topology induced by $\theta$ on $\S_0$
 but we do not have this yet.

\thrm Theorem 1.2."Let ${\frak f}\subset(\cl F.^*)_+$ be a folium.
Then it is the set of normal positive forms of the von Neumann 
algebra $\pi(\cl F.)''$ where $\pi=\mathop{\oplus}\limits_{\varphi
\in{\frak f}}\pi_\varphi.$  Conversely, the set of positive
normal forms of any von Neumann algebra is a folium."
See
 Haag, Kadison, Kastler in  [HKK].

This is quite useful, in that any folium can now be analyzed
as the normal state space of some concrete C*--algebra.

\item{\bf Remark} {\bf 1.3.}
For later use, we need to know about projections
associated with faces and folia.
Start with a pair  $(\cl F.,\,\S_0),$ where by Theorem~1.1
we now assume that $\S_0$  is a norm--closed face.
Corresponding to this, we know from Pedersen 3.6.11 [Pe]
that there is a projection $P\in\cl F.''$ which we now show 
how to construct. (To use Pedersen 3.6.11, we need to know that
a norm closed face $\S_0$ generates a cone $\r_+\S_0$
which is norm--closed and hereditary and in $(\cl F.'')_*,$ but it is quite
straightforward to verify this).
First define
$$\Gamma:=\set\varphi\in(\cl F.'')_*,|\varphi|\in\r_+\S_0.,$$
and this is in fact a left invariant vector space, by the
proof in  Pedersen 3.6.11. Then its annihilator $\Gamma^\perp\subset
\cl F.''$ is a $\sigma\hbox{--weakly}$ closed left
ideal, hence ${\Gamma^\perp\cap(\Gamma^\perp)^*}$ is a 
weak--operator closed hereditary subalgebra of  $\cl F.''.$
If  we denote its unit (which is a projection in  $\cl F.'')$
by $Q,$ then the desired projection we want is $P=\un-Q.$
To recover $\r_+\S_0$ from $P,$ we just take the set of
$\varphi\in(\cl F.''_*)_+=(\cl F.^*)_+$ such that
 $\varphi(P)=1.$
\chop
In the case that $\S_0$ is a folium (i.e. also invariant),
we find that $\Gamma$ is a two--sided invariant space, 
hence $\Gamma^\perp$ is a two--sided  $\sigma\hbox{--weakly}$ closed
ideal. It then follows from Pedersen~2.5.4~[Pe] that its unit
$Q\in\cl F.'\cap\cl F.'',$ and hence $P\in\cl F.'\cap\cl F.''.$
\chop
An obvious method by which one may think one can construct a host algebra, is
to take the algebra $\wt{\cl L.}:=C^*({P\cl F.P})\subset\cl F.''.$
Whilst this is certainly hereditary in $\cl F.'',$ and the states 
which uniquely extend from $\wt{\cl L.}$ to
$\cl E.:=C^*(\wt{\cl L.}\cup\cl F.)\subset\cl F.''$
will satisfy $\omega(P)=1,$ this is not enough to guarantee that
their restrictions to $\cl F.$ will be in $\S_0.$
This is because given a $\varphi\in\S(\cl F.),$
one can only conclude that $\varphi\in\S_0$
if its {\it normal} extension to $\cl F.''$ satisfies
$\wt\varphi(P)=1,$ and whilst for an $\omega\in\S(\cl E.)$
which extended from one on $\wt{\cl L.}$
we have $\omega(P)=1,$ we do not know 
that $\omega$ is the normal extension of its restriction
$\omega\restriction\cl F..$ Thus $\theta$ may not map
onto $\S_0$ for this choice $\wt{\cl L.}.$
By the previous remark we know there are folia without hosts,
so that we know the above procedure must sometimes fail.

\noindent
From Takesaki Prop.~2.17 (p129) [Tak] we know that if $\cl I.$
is a closed two--sided ideal of $\cl A.,$ then
$\pi(\cl I.)''=\pi(\cl A.)''$ for any representation
$\pi$ which is nondegenerate on $\cl I..$
This implies that if we have an ideal host $\cl L.\subset\cl E.$
for $(\cl F.,\,\S_0),$ then $\pi(\cl L.)''=\pi(\cl E.)''$
for the representation $\pi=\mathop{\bigoplus}\limits_{\omega\in\S(\cl L.)}
\pi_\omega.$ This fact leads us to suspect that $\pi(\cl L.)''=\pi(\cl F.)''$
and this is what we now want to prove, but we need a
  lemma first. We use Pedersen's notation $[\cdot]$ for
``closed linear span.'' 

\thrm Lemma 1.4." Let $(\cl F.,\,\S_0)$ have an ideal host
$\cl L.\subset\cl E.,$ and let $\pi:\cl E.\to\cl B.(\cl H.)$ be a representation.
Let $\cl H._e$ be the essential
subspace of $\pi(\cl L.).$
 Then for any vector  $\Omega\in\cl H._e$ we have
$$\big[\pi(\cl L.)\Omega\big]=\big[\pi(\cl F.)\Omega\big]\,.$$"
Denote $\cl H._\Omega:=\big[\pi(\cl L.)\Omega\big]$ and
$\pi_\Omega:=\pi\restriction\cl H._\Omega,$
then due to the fact that $\cl L.$ is an ideal host,
 $\pi_\Omega$ extends uniquely on the same space $\cl H._\Omega$
to $\f.$ Hence $\big[\pi(\cl F.)\Omega\big]=\big[\pi_\Omega(\cl F.)\Omega\big]
\subset\cl H._\Omega=\big[\pi(\cl L.)\Omega\big].$
We prove the reverse inclusion by contradiction.
Assume  it is not true, then there exists some
nonzero $\psi\in[\pi(\cl L.)\Omega]\subset\cl H._e$ such that
 $\psi\perp[\pi(\cl F.)\Omega].$
We have
$$\pi(\cl F.)\psi\perp\pi(\cl F.)\Omega
\eqno{(5)}  $$
because $\big(\pi(A)\psi,\,\pi(B)\Omega\big)
=\big(\psi,\,\pi(A^*B)\Omega\big)=0$ for all $A,\; B\in\cl F.$.\chop
Normalise: $\|\psi\|=1=\|\Omega\|$ and choose $\alpha,\,\beta\in\C$ such that
${|\alpha|^2+|\beta|^2}=1$ and $\alpha\not=0\not=\beta$ and define:
$$\eqalignno{\varphi:=&\alpha\Omega+\beta\psi\,,\qquad\omega_\varphi(A):=(\varphi,
\,\pi(A)\varphi)   \cr
\omega_{\psi}(A):=&(\psi,\,\pi(A)\psi)\,,
\qquad\omega(A):=(\Omega,\,\pi(A)\Omega)\;.\cr}$$
Now since $\varphi\in\cl H._e,$ $\omega_\varphi$ is nondegenerate on
$\cl L.,$ hence $\theta(\omega_\varphi\restriction\cl L.)=\omega_\varphi
\restriction\cl F.$ and so for all $A\in\cl F.$ we have:
$$\eqalignno{
\theta(\omega_\varphi)(A)&=\omega_\varphi(A)=
\Big(\alpha\Omega+\beta\psi,\,\pi(A)\big(\alpha\Omega+\beta\psi\big)
\Big)  \cr
&=|\alpha|^2\omega(A)+|\beta|^2\omega_{\psi}(A)
=\lambda\omega(A)+(1-\lambda)\omega_{\psi}(A)
\cr}$$
where we made use of the orthogonality~(5)
and we have set $\lambda:=|\alpha|^2\in(0,\, 1).$ Thus
$$\theta(\omega_\varphi)=\theta\big(\lambda\omega+(1-\lambda)\omega_{\psi}\big)\;.
\eqno{(6)}$$
However, for an element $L\in\cl L.$ we cannot use the orthogonality~(5),
and so we get
$$
\omega_\varphi(L)=\lambda\omega(L)+(1-\lambda)\omega_{\psi}(L)
+\bar\alpha\beta(\Omega,\,\pi(L)\psi)
+\alpha\bar\beta(\psi,\,\pi(L)\Omega)
\eqno{(7)}$$
We show that the last two terms can always be made nonzero by some choice of $L.$
If this were not the case, they must be zero for all $L=\gamma R$ where
$\gamma\in\C,$ $R=R^*\in\cl L..$ Then we have for the last two terms
of Equation~(7) that
$$ 2{\rm Re}\big[\alpha\bar\beta\gamma(\psi,\,\pi(R)\Omega)\big]=0 $$
for all $\gamma$ and $R,$ i.e. $(\psi,\,\pi(R)\Omega)=0$ for all
$R=R^*\in\cl L..$ However, $\cl L.$ is spanned by its selfadjoint elements,
thus $\psi\perp\pi(\cl L.)\Omega,$ and so
$$\psi\perp\left[\pi(\cl L.)\Omega\right]
\ni\psi$$
and thus $\psi=0$ which is a contradiction with our initial
assumption. 
Thus the last two terms of Equation~(7) are nonzero 
for some $L,$ i.e. on $\cl L.$ we have
$$\omega_\varphi\not=\lambda\omega+(1-\lambda)\omega_{\psi}\;.$$
This, together with Equation~(6) contradicts the assumption
 that $\theta$ is injective on  $\S(\cl L.)\,.$
Thus our initial assumption is wrong, so\chop
$\big[\pi(\cl L.)\Omega\big]\subseteq\big[\pi(\cl F.)\Omega\big]\,,$
and in fact we have equality
$\big[\pi(\cl L.)\Omega\big]=\big[\pi(\cl F.)\Omega\big]\,.$

\item{\bf Remarks.} {\bf (1)}
 Until now, we have used the standard notation
$\cl F.''$ for the universal von Neumann algebra of $\cl F.$
(not unique for concrete C*--algebras).
To avoid confusion in subsequent arguments, we will sometimes
explicitly indicate the universal representations, 
and our notation is that $\pi\s\f.:\f\to\cl B.(\cl H._\f)$
is the universal representation of $\f,$ i.e.
$\pi\s\f.=\mathop{\bigoplus}\limits_{\omega\in\S(\f)}\pi_\omega$
and $\f''\equiv\pi\s\f.(\f)''.$ Note that when $\cl L.$ is an ideal host
for $(\cl F.,\,\S_0),$ then $\pi\s\l.$ extends uniquely on
the same space to a representation of $\cl F.,$ and this implies
that $\pi\s\l.$ of $\f$ is a subrepresentation of $\pi\s\f..$
\item{\bf (2)}
One may try to generalise this lemma away from ideal hosts to hosts,
in which case we suspect that
 for any vector  $\Omega\in\cl H._e$ we have
$$\big[\pi(\cl L.)\Omega\big]\subseteq\big[\pi(\cl F.)\Omega\big]\,.$$
However, the proof so far eludes us. If one starts as in the proof by
assuming some nonzero $\psi\in{[\pi(\cl L.)\Omega]\Big\backslash
[\pi(\cl F.)\Omega]},$ then $\psi=\psi_0+\psi_1$ where
$\psi_0\in[\pi(\cl F.)\Omega]\perp\psi_1$ but we may have that
$\psi_i\not\in\cl H._e,$ even though $\psi\in\cl H._e,$ and
this causes problems.

Now we are ready to prove:
\thrm Theorem 1.5." Let $(\cl F.,\,\S_0)$ have an ideal host
$\cl L.,$ 
then $\pi\s\l.(\cl F.)''=\pi\s\l.(\cl L.)''\; (\equiv\cl L.'')$
where we use the same symbol   $\pi\s\l.$ for the unique extension
of it from $\cl L.$ to $\cl L.''\supset M(\cl L.)\supset\f.$"
In $\cl B.(\cl H._\l)$ let $A\in\pi\s\l.(\cl L.)'$
and let $B\in\cl F.,$ and recall that
$\pi\s\l.(B)\psi=\lim\limits_\alpha\pi\s\l.(BE_\alpha)\psi$
for all $\psi\in\cl H._{\l}$ and any approximate identity
$\{E_\alpha\}$ of $\cl L..$ Thus for all $\psi\in\cl H._\l$ we have
$$\big[A,\,\pi\s\l.(B)\big]\psi = \lim_\alpha
\big[A,\,\pi\s\l.(BE_\alpha)\big]\psi=0$$
because $A\in\pi\s\l.(\cl L.)'$ and $\cl L.$ is an ideal.
This is true for all $B\in\cl F.,$ and so $\pi\s\l.(\cl L.)'
\subseteq\pi\s\l.(\cl F.)'.$\chop
We now prove the reverse inclusion, and we do it by contradiction.
Assume that $\pi\s\l.(\cl F.)'\not=\pi\s\l.(\cl L.)'.$ Then since von Neumann algebras
are spanned by their projections, we can find a nontrivial projection
$P\in\pi\s\l.(\cl F.)'\backslash\pi\s\l.(\cl L.)'$ (otherwise, if all the projections of 
$\pi\s\l.(\cl F.)'$ were in $\pi\s\l.(\cl L.)'$ the algebras would be equal).
Recall that
$\pi\s\l.=\bigoplus\limits_{\omega\in\S(\cl L.)}\pi_\omega$,
so there must be some state in $\S(\cl L.)$, say $\omega_0$, such that
$$[\pi\s\l.(\cl L.),\, P]\,P_{\omega_0}\not=0=[\pi\s\l.(\cl F.),\, P]
\eqno{(8)}$$
where $P_{\omega_0}$ denotes the projection onto the subspace
$\cl H._{\omega_0}\subset\cl H._\l$ of the subrepresentation
$\pi_{\omega_0}:\cl L.\to\cl B.(\cl H._{\omega_0}).$
Let $\Omega_{\omega_0}$ be the normalised cyclic vector
for this representation. We claim that
$$
\|P\Omega_{\omega_0}\|\not=0\not=\|(\un-P)\Omega_{\omega_0}\|\;.$$
If $P\Omega_{\omega_0}=0,$ then
$$
0=\big[\pi\s\l.(\cl F.)P\Omega_{\omega_0}\big]=
\big[P\pi\s\l.(\cl F.)\Omega_{\omega_0}\big]=
\big[P\pi_{\omega_0}(\cl F.)\Omega_{\omega_0}\big]=P\cl  H._{\omega_0}$$
where we made use of Lemma~1.4, that $\Omega_{\omega_0}$ is 
cyclic for $\pi_{\omega_0}(\cl F.).$
But if $P$ annihilates $\cl  H._{\omega_0},$ it must commute with
$\pi\s\l.(\cl L.)$ on $\cl  H._{\omega_0},$ and this contradicts
Equation~(8), hence $P\Omega_{\omega_0}\not=0.$
\chop
Similarly, if $(\un-P)\Omega_{\omega_0}=0,$ then
by the same argument ${(\un-P)}$ commutes with
$\pi\s\l.(\cl L.)$ on $\cl  H._{\omega_0},$ hence so does $P.$
So, also ${(\un-P)\Omega_{\omega_0}}\not=0.$
Thus we can write:
$$\eqalignno{
\Omega_{\omega_0}&=P\Omega_{\omega_0}+(\un-P)\Omega_{\omega_0}  \cr
&=\alpha\cdot\Omega_P+\beta\cdot\Omega_{P^\perp}\qquad\hbox{where:} \cr
\Omega_P&:={P\Omega_{\omega_0}\over\|P\Omega_{\omega_0}\|}\;,\qquad
\Omega_{P^\perp}:={(\un-P)\Omega_{\omega_0}\over
\|(\un-P)\Omega_{\omega_0}\|}
\cr}$$
and where ${|\alpha|^2+|\beta|^2}=1$ and $\alpha\not=0\not=\beta.$
Adapting now the proof in Lemma~1.4, since we have that
${[\pi\s\l.(\cl F.),\, P]}=0$ we find for all $A\in\cl F.$:
$$\eqalignno{
\omega_0(A)&=\lambda\omega_P(A)+(1-\lambda)\omega_{P^\perp}(A)\qquad\hbox{where:} \cr
\omega_P(N)&=(\Omega_P,\,\pi\s\l.(N)\Omega_P)\;,\quad\hbox{and}\quad  \cr
\omega_{P^\perp}(N)&=(\Omega_{P^\perp},\,\pi\s\l.(N)\Omega_{P^\perp})\;,
\cr}$$
and we set $\lambda:=|\alpha|^2.$ That is, we have
$$\theta(\omega_0)=\theta\big(\lambda\omega_P+(1-\lambda)\omega_{P^\perp}\big)\;.
\eqno{(9)}$$
Now for an $L\in \cl L.$ we have similar to before:
$$\eqalignno{
\omega_0(L)=&
\Big(\alpha\Omega_P+\beta\Omega_{P^\perp},\,\pi\s\l.(L)\big(\alpha\Omega_P
+\beta\Omega_{P^\perp}\big)\Big) \cr
=&
\lambda\omega_P(L)+(1-\lambda)\omega_{P^\perp}(L)  \cr
&\qquad\quad+\bar\alpha\beta(\Omega_P,\,\pi\s\l.(L)\Omega_{P^\perp})
+\alpha\bar\beta(\Omega_{P^\perp},\,\pi\s\l.(L)\Omega_P)\;.&(10)\cr}$$
We prove that we can always find an $L\in\cl L.$ to make the last line nonzero.
If it is always zero, it is zero for all $L=\gamma R$ where $\gamma\in\C$
and $R=R^*\in\cl L..$ Thus
$$ 2{\rm Re}\big[\alpha\bar\beta\gamma(\Omega_{P^\perp},\,\pi\s\l.(R)\Omega_P)\big]=0
\qquad \forall\;\gamma\in\C,\;\;R=R^*\in\cl L.  $$
and thus $(\Omega_{P^\perp},\,\pi\s\l.(R)\Omega_P)=0$ for all $R=R^*\in\cl L..$
But $\cl L.$ is spanned by its selfadjoint elements, and so 
$$(\Omega_{P^\perp},\,\pi\s\l.(\cl L.)\Omega_P)=0\qquad\hbox{i.e.}\qquad
\pi\s\l.(\cl L.)\Omega_P\perp\pi\s\l.(\cl L.)\Omega_{P^\perp}\;.$$
Thus $\pi\s\l.(\cl L.)$ restricted to $[\pi\s\l.(\cl L.)\Omega_P]\oplus[\pi\s\l.(\cl L.)
\Omega_{P^\perp}]$
decomposes into two cyclic representations $\pi_P\oplus\pi_{P^\perp}.$
Since $\Omega_P$ is cyclic for $\pi_P(\cl L.),$ by Lemma~1.4
it is also cyclic for $\pi_P(\cl F.).$ Thus,
$$[\pi\s\l.(\cl L.)\Omega_P]=[\pi\s\l.(\cl F.)\Omega_P]
=[\pi\s\l.(\cl F.)P\Omega]=[P\pi\s\l.(\cl F.)\Omega]=P\cl H._{\omega_0}$$
and likewise we get $[\pi\s\l.(\cl L.)\Omega_{P^\perp}]=(\un-P)\cl H._{\omega_0}.$ So
for all $L\in\cl L.$ and $A\in\cl F.$ we have
$$\eqalignno{
P\pi\s\l.(L)\pi\s\l.(A)\Omega&=P\pi\s\l.(LA)(\alpha\Omega_P+\beta\Omega_{P^\perp})
=\alpha\pi\s\l.(LA)\Omega_P  \cr
&= \pi\s\l.(LA)P\Omega=\pi\s\l.(L)P\pi\s\l.(A)\Omega\;.  \cr}$$
Thus $[P,\,\pi\s\l.(L)]P_{\omega_0}=0$ for all $L\in\cl L..$
But this contradicts Equation~(8), hence the last line in
Equation~(10) is nonzero for some $L\in\cl L..$ Thus
$$\omega_0\not=\lambda\omega_P+(1-\lambda)\omega_{P^\perp}$$
and this, together with Equation~(9) now
contradicts the assumption that $\theta$ is injective on
$\S(\cl L.).$
Thus, the initial assumption was wrong, and we conclude
 $\pi\s\l.(\cl L.)'=
\pi\s\l.(\cl F.)',$ and hence  $\pi\s\l.(\cl L.)''=\pi\s\l.(\cl F.)''.$

The unique extension of $\pi\s\l.$ from $\cl L.$ to
$\cl F.$ is of course
just the representation $\pi\s\S_0.:=\bigoplus\limits_{\omega\in\S_0}\pi_\omega
\in{\rm Rep}(\cl F.),$ using  the definition of an ideal host.
\thrm Corollary 1.6."Let $\cl N.$ be a von Neumann algebra, and
let $\S_0$ be its set of normal states. If $\cl L.$ is an ideal host 
for the pair $(\cl N.,\,\S_0),$ then $\pi\s\S_0.(\cl N.)=\cl L.''$
and hence $\pi\s\S_0.(\cl N.)$ contains  $\cl L.$ as an ideal."
Recall the embedding $\cl N.\subset M(\cl L.)\subset\cl L.''.$
Note that if $\omega\in\S(\cl L.),$ then the unique extension
of $\pi_\omega$ to $\cl N.$ is (unitarily equivalent to)
$\pi_{\theta(\omega)},$ and this we see from
$$\theta(\omega)(N)=\lim_\alpha\omega(NE_\alpha)
=\lim_\alpha(\Omega_\omega,\,\pi_\omega(NE_\alpha)\Omega_\omega)
=(\Omega_\omega,\,\wt\pi_\omega(N)\Omega_\omega).$$
Thus, recalling that the universal representation of $\cl L.$
is $\pi\s\l.=\mathop{\bigoplus}\limits_{\omega\in\S(\cl L.)}\pi_\omega$
and that $\theta(\S(\cl L.))=\S_0,$ we conclude that the
unique extension of $\pi\s\l.$ to $\cl N.$ is
$\pi\s\l.\restriction\cl N.=\mathop{\bigoplus}\limits_{\omega\in\S_0}\pi_\omega
=\pi\s\S_0..$
Thus $\pi\s\S_0.$ on $\cl N.$ is a normal representation, and so
$$\pi\s\S_0.(\cl N.)=\pi\s\S_0.(\cl N.'')=\pi\s\S_0.(\cl N.)''=\cl L.''$$
where we used Theorem~1.5 for the last equality.

\thrm Corollary 1.7." If  $\cl L.$ is an ideal host for a pair
$(\f,\S_0),$ then $\theta:\S(\l)\to\S_0$ maps the pure states of 
$\l$ to pure states on $\f.$
Hence all extreme points of $\S_0$ are pure."
By Theorem~1.5 we have $\pi_{\S_0}(\f)''=\pi_{\S_0}(\l)''.$
Let $\omega_0\in\S_0,$ then since $\pi_{\S_0}=\bigoplus\limits_{\omega\in\S_0}
\pi_\omega,$ we have that the restriction map
$R:\pi_{\S_0}(\f)\to\pi_{\omega_0}(\f)$ by $R(\pi_{\S_0}(F)):=\pi_{\omega_0}(F)$
for all $F\in\f$ is a normal *--homomorphism, hence it extends to
$\pi_{\S_0}(\f)''=\pi_\l(\l)''$ and so
$$\eqalignno{R\big(\pi_{\S_0}(\f)''\big)&=R\big(\pi_{\S_0}(\f)\big)''
=\pi_{\omega_0}(\f)''\cr
&=R\big(\pi_\l(\l)\big)''=\pi_{\omega_0}(\l)''\,.\cr}$$
Or to be more notationally precise, 
$\pi_{\omega_0}(\f)''=\pi_{\theta^{-1}(\omega_0)}(\l)''.$
Now $\theta$ is an affine bijection, so it restricts to a
bijection between the pure states of $\S(\l)$ and the extreme
points of $\S_0.$
If $\varphi$ is a pure state on $\l,$ then
$\b(\cl H._{\varphi})=\pi_{\varphi}(\l)''=\pi_{\theta(\varphi)}(\f)'',$
i.e. $\theta(\varphi)$ is also pure on $\f.$

This last corollary now severely limits the class of folia for which
ideal hosts exist, indeed, we can quickly prove many obstruction theorems,
e.g. the next one.
\thrm Corollary 1.8." If  $\cl N.$ is a simple factor and
$\S_0$ is  its folium of normal states, then there is no ideal host
for the pair $(\cl N.,\,\S_0).$"
it suffices by Corollary~1.7 to observe that $\S_0$ contains no
pure states. For if an $\omega\in\S_0$ were pure, then
using the fact that $\pi_\omega$ is normal, we see
that $\pi_\omega(\cl N.)=\pi_\omega(\cl N.)''=\b(\cl H._\omega).$
Since $\cl N.$ is simple, $\pi_\omega$ is an isomorphism
hence $\cl N.\cong\b(\cl H._\omega)$ and the latter is not simple.
This is a contradiction, so $\S_0$ has no pure states.

\item{\bf Remarks.}{\bf (1)}
Thus, since we know from Kadison and Ringrose 6.8.4 and 6.6.5 [KR],
that every finite factor and each countably decomposable type III factor
is simple, Corollary~1.8 shows that there are many von Neumann 
algebras for which there is no ideal host for its normal states.
From Corollary~1.6 we see that we
can only expect ideal hosts for a von Neumann algebra
which has a norm closed proper  ideal which is weak operator dense.
If $\cl F.$ is a C*--algebra, but not a von Neumann algebra, 
 Theorem~1.5 just tells us
that we should be looking for an ideal host $\cl L.$ in
the von Neumann algebra $\pi\s\S_0.(\cl F.)'',$ but by  Corollary~1.7
we also know that this may fail, unless all the extreme points
of $\S_0$ are pure.
\item{\bf (2)}
 Theorem~1.5 states a ``weak'' uniqueness,
in that it claims that all ideal hosts for the same pair
must  have the same universal algebra
$$\l''=  \pi_{\S_0}(\f)''.$$
Theorem 1.5 also tells us where all ideal hosts reside, viz
$\cl L.\subset\pi\s\S_0.(\cl F.)''\subset\cl B.(\cl H._{\S_0}),$
and we can make this even more precise:
Since $\l$ is an ideal host, $\S_0$ is a folium, and so 
for the projection associated with $\S_0$ we have
$P\in\f'\cap\f''$ (cf. Remark~1.3)
and  $\S_0=\set\omega\in\S(\f),\wt\omega(P)=1.$
where $\wt\omega$ denotes the normal extension of $\omega$
from $\f$ to $\f''.$
Thus in $\pi\s\f.$ we have $P\Omega_\omega=\Omega_\omega$
for all $\omega\in\S_0,$ hence by $P\in\f'\cap\f''$
we see $P F\Omega_\omega= F\Omega_\omega$
for all $F\in\f''.$ Thus $P$ is the projector onto
$\bigoplus\limits_{\omega\in\S_0}\cl H._\omega=\cl H._{\S_0}
\subset\cl H._\f,$ hence $P\pi\s\f.(\cl F.)''=
\pi\s{\S_0}.(\cl F.)''$ and so we conclude
from $P\in\f''$ that
$$\l\subset\pi\s\S_0.(\cl F.)''\subset\pi\s\f.(\cl F.)''=\f''\,.
\eqno{(11)}$$
\item{\bf (3)}
Since for an ideal host $\l$ for a pair $(\f,\,\S_0)$
we know by Theorem~1.5 that
$\xi(\cl F.)\subseteq M(\cl L.)\subset\l''=\pi\s\l.(\cl F.)'',$
we conclude that the embedding $\xi:\f\to M(\l)$ is precisely
$\pi\s\S_0.\restriction\f.$
So $\xi$ is faithful iff $\pi\s\S_0.\restriction\f$ is
faithful iff for each $F\in\f$ we have $\omega(F)\not=0$
for some $\omega\in\S_0.$
\item{\bf (4)}
We can easily adapt the proofs of 1.4 and 1.5
to prove a 
Stone--Weierstrass theorem for von Neumann algebras, i.e.
if a von Neumann algebra contains a sub--von Neumann algebra which
separates its normal states, then they are equal. 
However, this fact has also a very short proof via the bipolar
theorem (Private communication with R.Longo).
\item{\bf (5)}
By Corollary 1.7 we can also find topological goups which have no group
algebras (i.e. ideal hosts for the the pair ${\left(C^*(G_d), \S_0\right)}$
where $\S_0$ are states $\omega$ for which $g\to\omega(\delta_g)$
is continuous). For example, let $G=L^\infty(\r)$ with the group
operation being addition, and with the strong operator topology 
w.r.t. its representation as multiplication operators on $L^2(\r).$
Then there is a topological isomorphism $L^\infty(\r)\cong C(\r_s)$
where $\r_s$ is $\r$ compactified and
endowed with a suitable hyperstonean topology
 (cf. proof of Theorem~III.1.18 [Tak], or Theorem~2.1
below). Any irreducible representation 
of $G$ must be a character, hence point evaluation
on $C(\r_s),$ and this cannot be continuous w.r.t. 
the strong operator topology, because points are still
of measure zero w.r.t. the extension of the Lebesgue
measure to $\r_s\,.$ Thus the folium $\S_0$ of states of
$C^*(G_d)$ which are continuous w.r.t. the topology of $G$ 
contains no pure states, hence by Corr.~1.7 we
conclude that G has no group algebra.

\noindent

To conclude this section, we would like to make precise the relation between the 
representations of a host $\l$ for a pair $(\f,\S_0)$ and the 
representations of $\f.$ Denote the normal representations of
$\pi\s\S_0.(\f)$ by ${\rm Rep}\s\S_0.\f.$ By Theorem~1.5 we know
that $\l''\subseteq\pi\s\S_0.(\f)'',$ hence for each
$\pi\in{\rm Rep}\s\S_0.\f$ we can construct a representation
$\Lambda(\pi)\in{\rm Rep}\,\l$ by first extending $\pi$ via
strong operator continuity to a representation
$\wt\pi\in{\rm Rep}\big(\pi\s\S_0.(\f)''\big),$ and then
defining $\Lambda(\pi)$ as $\wt\pi\restriction\l$ restricted to its
essential subspace. This produces a map
$\Lambda:{\rm Rep}\s\S_0.\f\to{\rm Rep}\,\l.$
\thrm Theorem 1.9." If $\l$ is an ideal host for the pair
$(\f,\S_0),$ then the map
$\Lambda:{\rm Rep}\s\S_0.\f\to{\rm Rep}\,\l$
is a bijection which takes irreducible representations to
irreducible representations. In the case that $\l$ is merely
a host, $\Lambda$ is a bijection modulo unitary equivalence.
Its inverse is via inducing of representations;-
$\Lambda^{-1}(\{\pi\})=\big\{{\rm Ind}_\l^\f(\pi)\big\}$
where $\{\cdot\}$ denotes unitary equivalence classes, 
and the induction is done via the 
right\ab$\l\hbox{--rigged}$left\ab$\f\hbox{--module}$
$\cl M.:=[\f\l]\subset\cl E.$ with rigging map
${\langle u,\,v\rangle}:=u^*v\in\l$ for all
$u,\, v\in\cl M..$"
If $\l$ is an ideal host, $\l''=\pi\s\S_0.(\f)'',$
and now as both $\l$ and $\pi\s\S_0.(\f)$ are strong operator dense
in $\l'',$ it is obvious that each uniquely determines a normal
representation, and so the proof for this case ends here.\chop
For the case of $\l$ just a host, let $(\pi,\cl H.)\in
{\rm Rep}\s\S_0.\f,$ so $\Lambda(\pi)$ is $\wt\pi\restriction\l$
restricted to its essential subspace $[\wt\pi(\l)\cl H.].$
Now for our proof, we will first construct 
${\rm Ind}_\l^\f(\Lambda(\pi)),$ show it is unitary equivalent to
$\pi,$ and then show that every $\pi\in{\rm Rep}\s\S_0.\f$
is unitarily equivalent to some 
${\rm Ind}_\l^\f(\gamma),$ $\gamma\in{\rm Rep}\,\l.$
Following Fell and Doran~XI.4.12~[FD] or Rieffel~[Ri],
consider the
right\ab$\l\hbox{--rigged}$left\ab$\f\hbox{--module}$
$$\cl M.:=[\f\l]\subset\cl E.,\qquad
{\langle u,\,v\rangle}:=u^*v\in\l\quad\forall\;
u,\, v\in\cl M.$$
where we used the fact that $\l$ is hereditary in $\cl E.$
to conclude $u^*v\in\l.$ Since $\l$ is a C*--algebra, every
representation of it is inducible via $\cl M.$ (cf.
XI.4.12~[FD]). Now construct
$\rho={\rm Ind}_\l^\f(\gamma)\in{\rm Rep}\,(\f)$ as follows.
On $\cl M.\otimes\cl H.$ define a pre--inner product
$(\cdot,\cdot)_0$ by
$$\big(s\otimes\xi,\, t\otimes\eta\big)_0:=
\left(\gamma\big(\langle t,s\rangle\big)\xi,\,\eta\right)=
\left(\gamma\big(t^*s\big)\xi,\,\eta\right)\eqno{(12)}$$
and define from $\cl M.\otimes\cl H.$ the Hilbert space
$$\cl K.:=\overline{\cl M.\otimes\cl H.\big/\ker(\cdot,\cdot)_0}$$
where closure is obviously w.r.t. $(\cdot,\cdot)_0.$
\def\ttensor{\mathrel{\wt{\mathord{\otimes}}}}
Denote the image of an elementary tensor $s\otimes\xi$
in $\cl K.$ by $s\ttensor\xi,$ and define the representation
$\rho:\f\to\cl B.(\cl K.)$ by
$$\rho(F)(s\ttensor\xi):=Fs\ttensor\xi\qquad\forall\;
s\in\cl M.,\;\xi\in\cl H.,\; F\in\f\;.$$
If we let $\gamma=\Lambda(\pi),$ then Equation~(12) becomes
$$\big(s\otimes\xi,\, t\otimes\eta\big)_0:=
\left(\wt\pi(s)\xi,\,\wt\pi(t)\eta\right)$$
and so we can identify $\cl K.$ with the subspace
${[\wt\pi(\cl M.)\cl H.]}={[\wt\pi(\f\l)\cl H.]},$
via the unitary $U(s\otimes\xi):=\wt\pi(s)\xi.$
Moreover we have for all $s\ttensor\xi,$ $t\ttensor\eta$ that
$$\eqalignno{\left(\rho(F)(s\ttensor\xi),\,t\ttensor\eta\right)&=
\left(Fs\ttensor\xi,\,t\ttensor\eta\right)=
\left(\wt\pi(Fs)\xi,\,\wt\pi(t)\eta\right)  \cr
&=\left(\pi(F)\cdot\wt\pi(s)\xi,\,\wt\pi(t)\eta\right)\cr}$$
and so $(\rho,\,\cl K.)$ is unitarily equivalent to
$\pi\restriction\f$ on ${[\wt\pi(\cl M.)\cl H.]}.$
Now application of Lemma~1.10 (proven below) to
$\wt\pi\in{\rm Rep}\,\cl E.$ implies that we have 
$\cl H.={[\wt\pi(\cl M.)\cl H.]}$ iff $\pi\in
{\rm Rep}\s\S_0.\f,$ and the latter is what we assumed at the start.
Thus $\rho$ is unitarily equivalent to $\pi.$
Furthermore, we see above that the induction process produce 
representations $\wt\pi$ such that $\cl H.={[\wt\pi(\cl M.)\cl H.]},$
hence the image under induction via $\cl M.$ of ${\rm Rep}\,\l$
is ${\rm Rep}\s\S_0.\f,$ using Lemma 1.10 again.

\thrm Lemma 1.10."A representation $(\pi,\cl H.)\in{\rm Rep}\,\cl E.$ 
satisfies 
$\cl H.={[\pi(\cl M.)\cl H.]}$ iff ${\pi\restriction\f}$ is normal
with respect to $\pi\s\S_0.(\f).$"
Let $\cl H.={[\pi(\cl M.)\cl H.]}$ and  choose a normalised vector
$\psi\in\pi(\l)\cl H.,$ i.e. 
$\psi=\pi(L)\xi,$ $\|\psi\|=1.$
Let $\omega_\psi$ denote the associated vector state on $\cl E.,$
then since $\|{\omega_\psi\restriction\l}\|=1$ and $\l$ is hereditary,
it is everywhere determined on $\cl E.$ by its values on $\l.$
Thus ${\omega_\psi\restriction\f}=\theta({\omega_\psi\restriction\l})
\in\S_0\subset\big(\pi\s\S_0.(\f)\big)_*.$
Since the normal functionals of any representation is a folium,
hence invariant under conjugation, we conclude that also
 ${\omega_{\pi(F)\psi}\restriction\f}={\omega\s\pi(FL)\xi.\restriction\f}
\in\big(\pi\s\S_0.(\f)\big)_*$ for all $L\in\l,$ $F\in\f$ and
$\xi\in\cl H.,$ i.e. 
 ${\omega\s\pi(\f\l){\cl H.}.\restriction\f}
\in\big(\pi\s\S_0.(\f)\big)_*.$ 
Since by assumption $\pi(\f\l)\cl H.$ spans a dense subspace of $\cl H.,$
it follows from Kadison and Ringrose~7.1.15~[KR] that ${\pi\restriction\f}$
is normal with respect to $\pi\s\S_0.(\f).$\chop
Conversely, let ${\pi\restriction\f}$ be
 normal with respect to $\pi\s\S_0.(\f),$ and assume that
$\cl H.\not={[\pi(\cl M.)\cl H.]},$ i.e. there is some
$\psi\perp
{[\pi(\cl FL.)\cl H.]}.$ First, we show that 
 ${\omega_\psi\restriction\f}
\not\in\big(\pi\s\S_0.(\f)\big)_*.$ If not, then $\theta^{-1}(
{\omega_\psi\restriction\f})\in\S(\l),$ i.e. the normal extension
$\wt\omega_\psi$ to $\f''$ restricts to a state on $\l.$
Since this normal extension on $\cl E.$ is just
$\wt\omega_\psi(A)=(\psi,\pi(A)\psi),$ we conclude from
the given $\psi\perp[\pi(\f\l)\cl H.]\supset\pi(\l)\psi$
that $\wt\omega_\psi(\l)=0,$ which contradicts the fact that it must be
a state on $\l.$ Thus
 ${\omega_\psi\restriction\f}
\not\in\big(\pi\s\S_0.(\f)\big)_*$. 
Now we know by the first part that
 ${\omega\s\pi(\f\l){\cl H.}.\restriction\f}
\in\big(\pi\s\S_0.(\f)\big)_*$, and in fact since the normal functionals
is a norm--closed folium, we have
 ${\omega\s[\pi(\f\l){\cl H.}].\restriction\f}
\in\big(\pi\s\S_0.(\f)\big)_*$, so since
$\cl H.=\C\psi\oplus[\pi(\f\l)\cl H.]$ it is impossible to find a dense subspace
$\cl S.\subset\cl H.$ such that ${\omega_\varphi\restriction\f}
\in\big(\pi\s\S_0.(\f)\big)_*$ for all $\varphi\in\cl S.,$
hence by Kadison and Ringrose~7.1.15~[KR] ${\pi\restriction\f}$
cannot be normal w.r.t. $\pi\s\S_0.(\f).$
This contradicts our hypothesis, hence 
$\cl H.={[\pi(\cl M.)\cl H.]}.$


\beginsection 2. Ideal hosts up to a central algebra.

Above we saw that a pair $(\cl F.,\,\S_0)$ with $\S_0$ a folium
without pure states, has no ideal host.
A particularly bad case of this, is the pair 
${\big(L^\infty(X,\mu),\,\S_N\big)}$ where $\mu$ has no discrete
part, and where $\S_N$ denotes the set of normal states, i.e. the 
measures absolutely continuous w.r.t. $\mu.$ 
In this case $\S_N$ does not even have extreme points, because
for any measure $\nu$ absolutely continuous w.r.t. $\mu,$ 
we only need to subdivide its support to write it as
a convex combination of other probability measures in this class.
Only when ${\rm supp}(\nu)$ is a point can we not do this,
and this case does not occur since $\mu$ has no discrete part.
In this section we want to argue that this example is 
symptomatic of the general case, in that if a pair 
 $(\cl F.,\,\S_0)$ has no ideal host, it is because
${\big(L^\infty(X,\mu),\,\S_N\big)}$ is embedded in it,
 and it acts as an obstruction.
To be precise about the embedding, we will show that for any pair
 $(\cl F.,\,\S_0)$ with $\S_0$ a folium,
 we can always find a quasi--host $\wl$ in the sense of the next
definition:
\item{\bf Def.}
Given a pair
 $(\cl F.,\,\S_0)$ where $\S_0$ is a folium,
 a {\it quasi--host} for it, is a C*--algebra $\wl$
and two embeddings $\f\subset M(\wl)$ and $L^\infty(X,\,\mu)
\subset ZM(\wl)$ for some measure space $(X,\,\mu)$
such that ${\wt\S_\mu\restriction\f}=\S_0$ where
$\S_\mu:=\set\omega\in\S(\wl),\wt\omega\restriction L^\infty(X,\mu)
\;\hbox{is normal}.$ and moreover,\chop
${\wt\S_\mu\restriction C^*(\f\cup L^\infty(X,\mu))}$ defines
an injection for $\S_\mu.$

\noindent
We will show that the given pair
 $(\cl F.,\,\S_0)$ has no ideal host if the measure $\mu$ is
purely continuous. This is what we mean by saying
$L^\infty(X,\mu)$ acts as an obstruction.

\item{\bf Example.} Consider the von Neumann algebra
$\f={L^\infty(X,\,\mu)\mathrel{\overline{\mathord{\otimes}}}
\cl B.(\cl H.)}$ acting on the Hilbert space
${L^2(X,\,\mu)\otimes\cl H.},$ where 
${\overline{\mathord{\otimes}}}$ denotes the W*--tensor product
(cf. 11.2~[KR]). Its pure states consists of product states
$\omega_1\otimes\omega_2$ such that $\omega_i$ are both pure.
Thus, if we take the pair $(\f,\S_0)$ where $\S_0$ are the normal
states of $\cl F.,$ and assume that $\mu$ has no discrete part,
then $\S_0$ has no pure states, hence this pair has no ideal host.
Nevertheless, the algebra $\wt\l=
{L^\infty(X,\,\mu)\otimes
\cl K.(\cl H.)}$ is a quasi--host for $(\f,\,\S_0).$

\noindent
To start the analysis, let 
 $(\cl F.,\,\S_0)$ be a pair with $\S_0$ a folium, then by Theorem~1.2
we construct the representation $\pi_{\S_0}=\bigoplus\limits_{\omega
\in\S_0}\pi_\omega$ and identify $\S_0$ with the normal states
of the concrete algebra $\pi_{\S_0}(\f),$ hence with the normal states
 of the von Neumann
algebra $\pi_{\S_0}(\f)''.$
We first analyze a single cyclic component $\pi_\omega,$ $\omega
\in\S_0$ of the direct sum.
 We denote the von Neumann algebra
$\pi_\omega(\f)''$ by $\cl N.,$ and its set of normal states by
$\S_N.$

Central to the following constructions, is the usual decomposition theory
with respect to some commutative subalgebra $\cl C.\subset\cl N.'$ (cf.
Takesaki [Tak]), however, since the maps involved
are only defined up to $\mu\hbox{--negligible}$ sets on some 
measure space $(X,\,\mu),$ and we will actually need everywhere defined maps,
we now redo some of the basic constructions in order to remedy this.
\thrm Theorem 2.1."Let $\cl F.$ be a C*--algebra with a fixed state
$\omega\in\S(\cl F.),$ and a commutative
 unital C*--algebra $\cl C.\subset
\pi_\omega(\cl F.)',$ then the Gel'fand isomorphism
$\Phi:\cl C.''\to L^\infty(X,\,\mu)$ equips the spectrum $X$
of $\cl C.$ with a probability measure $\mu,$ $\supp\mu=X,$ and there is 
a map $\psi:X\to\S(\cl C.')$ such that\chop
(\rn1) $\;$ the map $x\to\psi_x(A)$ is in $L^\infty(X,\,\mu)$ for all
$A\in\cl C.',$ and if $\cl C.$ is a von Neumann algebra,
$x\to\psi_x(A)$ is in $C(X)$ for all $A\in\cl C.',$\chop
(\rn2) $\;$ $\wt\omega(A):=(\Omega_\omega,\, A\Omega_\omega)
=\int_X\psi_x(A)\,d\mu(x)$ $\forall\; A\in\cl C.',$\chop
(\rn3) $\;$ $\psi_x(C\cdot A)=\Phi(C)(x)\cdot\psi_x(A)$
$\forall\; C\in\cl C.''$, $A\in\cl C.'.$"
(This proof is based on Takesaki 6.23, p241 [Tak])\chop
Since $\Omega_\omega$ is cyclic for $\pi_\omega(\cl F.),$
it is separating for $\cl C.$ and $\cl C.'',$ so
$\wt\omega\rest\cl C.$ is a faithful state of $\cl C..$
Thus by the Riesz representation theorem there is a
Borel measure $\mu$ on the compact set $X$ such that
$$\wt\omega(C)=\int_X\Phi(C)(x)\, d\mu(x)\;,\qquad
C\in\cl C.$$
with $\Phi:\cl C.\to C(X)$ the Gel'fand isomorphism.
Moreover $\mu$ is a probability measure since $\wt\omega(\un)=1$
and $\supp\mu=X$ because $\wt\omega$ is faithful.
Since $\Omega_\omega$ separates $\cl C.,$ we can consistently
define a unitary $U:[\cl C.\Omega_\omega]\to L^2(X,\,\mu)$ by
$U(C\Omega_\omega):=\Phi(C)$ which produces the representation
$\wt\Phi:\cl C.\to\cl B.(L^2(X,\,\mu))$ by $\wt\Phi(C)
:=UCU^{-1}=T\s\Phi(C).$ in terms of multiplication operators
$\set T_f,f\in C(X)..$  Since unitary 
conjugation is a normal map, $\wt\Phi$ extends to 
$\cl C.''\subset\pi_\omega(\cl F.)'$ and $\wt\Phi(\cl C.)''
=\wt\Phi(\cl C.''),$ i.e.  $\wt\Phi(\cl C.'')=
\set T_f,f\in L^\infty(X,\,\mu)..$
If $\cl C.$ is already a von Neumann algebra, this simplifies to
$$\set T_f,f\in C(X).=\wt\Phi(\cl C.)=\wt\Phi(\cl C.)''
=\set T_f,f\in L^\infty(X,\,\mu).\,,$$
so for each $f\in L^\infty(X,\,\mu)$ we can find a $\wt f\in C(X)$
such that $T_f=T\s\wt f.$ on $L^2(X,\,\mu),$ producing the
isomorphism $C(X)\cong L^\infty(X,\,\mu)$ (this isomorphism
also occurs in the proof of Theorem~III.1.18~[Tak]). Henceforth we will
blur the distinction between $\Phi$ and $\wt\Phi$ and always make the
identification with $C(X)$ if $\cl C.$ is a von Neumann algebra.\chop
Let $e$ be the projection of $\cl H._\omega$ onto 
\def\oo{\Omega_\omega}
$[\cl C.\Omega_\omega] =[\cl C.''\oo]$, then $e\in\cl C.'$ because
$\cl C.$ preserves $[\cl C.\oo]$. 
Define $\Upsilon:\cl C.''\to\cl B.\big([\cl C.\oo]\big)$ by
$\Upsilon(C):=e\, C$, then $\Upsilon(\cl C.'')=e\, \cl C.''$
is maximally commutative in $\cl B.\big([\cl C.\oo]\big)$ because it
has a cyclic vector $\oo$, cf. Takesaki Corr. 1.3, p104 [Tak].
Moreover $\Upsilon$ is injective because $\oo$ is separating for $\cl C.''$.
Now $(e\cl C.)'\rest[\cl C.\oo]=e\, \cl C.'e\subset\cl B.([\cl C.\oo])$
(easily verified, but also in Takesaki 3.10 [Tak]),
hence $e\, \cl C.'e\subset e\, \cl C.''$ because the latter is maximally
commutative in $\cl B.([\cl C.\oo])$ (also by Takesaki 3.10).
Define $\delta:\cl C.'\to\Upsilon(\cl C.'')$ by 
$\delta(A):=eAe\in e\, \cl C.''$. Clearly $\delta(\un)=e$ and $\delta$ is
positive and extends $\Upsilon$.
Next define $\psi:X\to\S(\cl C.')$ by 
$$\psi_x(A):=\Phi(\Upsilon^{-1}(\delta(A)))(x)$$
and note $\Phi(\Upsilon^{-1}(\delta(A)))$ is in $L^\infty(X,\,\mu)$, so $\psi$
is defined $\mu\hbox{--almost}$ everywhere, and if $\cl C.$ is a 
von Neumann algebra we can identify $\Phi(\Upsilon^{-1}(\delta(A)))$
with an element of $C(X)$ in which case $\psi$ is everywhere
defined. By positivity and normalisation of the various maps each $\psi_x$
is a state. Now
$$\leqalignno{\wt\omega(\Upsilon^{-1}(\delta(A)))&=\int_X\Phi(\Upsilon^{-1}
(\delta(A)))(x)\; d\mu(x)\cr
&=\int_X\psi_x(A)\; d\mu(x)\cr
\wt\omega(A)=(\oo,\,\,A\oo)&=(\oo,\,eAe\oo)\qquad\hbox{using $e\oo=\oo$}\cr
&=(\oo,\,\delta(A)\oo)=(\oo,\,\Upsilon^{-1}(\delta(A))\oo)\cr
&=\wt\omega(\Upsilon^{-1}(\delta(A)))=\int_X\psi_x(A)\, d\mu(x)\;.\cr}$$
For (\rn2), observe that for $C\in\cl C.''$, $A\in\cl C.'$ we have
$$\psi_x(C\cdot A)=\Phi(\Upsilon^{-1}(\delta(CA)))(x)
=\Phi(C)(x)\cdot\psi_x(A)$$
simply using the fact that $\Phi,\,\Upsilon^{-1}$ and $\delta$ are homomorphisms
and by definition $\Upsilon^{-1}(\delta(C)=C$ for $C\in\cl C.''$.

\itemitem{\bf Remarks.}{\bf (1)} A very important point here is that for $\cl C.$
a von Neumann algebra, $\psi:X\to\S(\cl L.')$ is everywhere defined. This comes 
from the identification which $\Phi$ makes between $C(X)$ and $L^\infty(X,\,\mu).$
The spectrum $X$ of $\cl C.$ in this case is of course hyperstonean
hence extremely disconnected (cf. [Tak]).
\itemitem{\bf (2)} Note that no separability assumption was needed here.
\itemitem{\bf (3)} The map $\psi$ and measure $\mu$ produces precisely the
orthogonal barycentric decomposition of $\wt\omega$ used in decomposition theory
in e.g. Takesaki [Tak] or Bratteli and Robinson [BR].
One uses $\psi$ to write $\mu$ as a measure on the state space
$\S(\cl F.).$
\itemitem{\bf (4)} By restricting $\wt\omega$  and $\psi_x$
to any subalgebra of
$\cl C.'$ (e.g. $\pi_\omega(\cl F.)$ or $\pi_\omega(\cl F.)''$), we obtain
a decomposition theorem for states on these.
Note that the usual theory assumes that $\cl C.$ is a von Neumann algebra.
\medskip
\item{\bf Def.} Given the data $\cl F.,\;\omega,\;\cl C.$ of theorem 2.1
assume that $\cl C.$ is a von Neumann algebra
with subsequent map $\psi:X\to\S(\cl C.')$ and define the following two 
bundles on $X$:\chop
$\bullet$ $\cl H.(\psi)$ is the bundle with projection  $p:\cl H.(\psi)\to X$
by $p^{-1}(x)=\cl H.\s\psi_x.$,\chop
$\bullet$ $\cl B.(\psi)$ the bundle with $q:\cl B.(\psi)\to X$ by
$q^{-1}(x)=\cl B.(\cl H.\s\psi_x.)$.

\noindent
At this point there is no topology, so the total spaces of these bundles
are just the unions of their fibres. 
Clearly the sections $\Gamma(\cl B.(\psi))$ act pointwise on the sections
$\Gamma(\cl H.(\psi))$. 
Now there are some canonical families of sections:
\item{$\bullet$} $\Omega\in\Gamma(\cl H.(\psi))$ is the section
$\Omega(x):=\Omega\s\psi_x.$,
\item{$\bullet$} $\Pi(A)\in\Gamma(\cl B.(\psi))$ is the section
$\Pi(A)(x):=\pi\s\psi_x.(A)$, $A\in\cl C.'$.

\noindent Clearly the latter specialises on $\cl C.$ to
\chop $\bullet$ $\Pi(C)(x)\varphi
=\Phi(C)(x)\;\varphi\qquad\forall\;\varphi\in\cl H.\s\psi_x.,\; C\in\cl C.$.
\medskip\noindent
Let $c=\Pi(C)\Omega$, $d=\Pi(D)\Omega$ with $C,\; D\in\cl C.'$, then
$$\big(c(x),\; d(x)\big)\s{\cl H.}_{\psi_x}.=
\big(\Pi(C)\Omega(x),\; \Pi(D)\Omega(x)\big)\s{\cl H.}_{\psi_x}.=
\psi_x(C^*D)$$
which is a continuous function in $x$ by theorem 2.1. Obviously
 $\Pi(\cl C.')\Omega$ is a linear space, and we have
that $x\to
\big(c(x),\; d(x)\big)\s{\cl H.}_{\psi_x}.$ is continuous, hence 
integrable (since $X$ compact, $\mu$ a probability measure).
Thus we can equip $\Pi(\cl C.')\Omega$ with the inner product
$$(c,\, d):=\int_X
\big(c(x),\; d(x)\big)\s{\cl H.}_{\psi_x}.\, d\mu(x)\;.$$
and hence form the Hilbert space
$$\cl H._\Gamma:=\overline{\Pi(\cl C.')\Omega\big/\ker(\cdot,\cdot)}$$
with factorisation map $\kappa:\Pi(\cl C.')\Omega\to\cl H._\Gamma$.
\thrm Theorem 2.2." $\Pi(\cl C.')$ lifts through $\kappa$ to define
a representation $\pi:\cl C.'\to\cl B.(\cl H._\Gamma)$. Then 
$\kappa(\Omega)$ is cyclic for $\pi(\cl F.)$, and there is a unitary 
$U:\cl H._\omega\to\cl H._\Gamma$ which intertwines the representation
$\big(\pi\s\tilde\omega.,\,\cl H.\s\tilde\omega.,\,\Omega\s\tilde\omega.\big)$
of $\cl C.'$ with 
$\big(\pi,\,\cl H._\Gamma,\,\kappa(\Omega)\big)$."
Obviously $\Pi(\cl C.')$ preserves $\Pi(\cl C.')\Omega$;- we show that it
preserves $\ker(\cdot,\cdot)$. Let $\varphi=
\Pi(A)\Omega\in
\ker(\cdot,\cdot)$, i.e.
$$\eqalignno{0=\|\varphi\|^2=\int_X
\big(\Omega(x),\,\Pi(A^*A)(x)\,\Omega(x)\big)\; d\mu(x)
=\int_X
\psi_x(A^*A)\; d\mu(x)\cr}$$
Thus $\psi_x(A^*A)=0$ by positivity and continuity of this map.
Let $B\in\cl C.'$, then
$$\|\Pi(B)\varphi\|^2=\int_X\psi_x(A^*B^*BA)\, d\mu(x)\;.$$
By the Cauchy--Schwartz inequality $\psi_x(A^*B^*BA)=0$
using \chop $\psi_x(A^*A)=0$. Thus $\Pi(B)$ preserves
$\ker(\cdot,\cdot)$ and hence lifts through $\kappa$.\chop
To show that $\kappa(\Omega)$ is cyclic for $\Pi(\cl F.)$, assume
the converse, i.e. there is a sequence $\{\varphi_n\}\subset
\Pi(\cl C.')\Omega$ converging with respect to the seminorm
$\|\varphi\|:=\left[\int_X\|\varphi(x)\|\s{\cl H._{\psi_x}}.d\mu(x)\right]^{
1/2}$
and $\big(\varphi_n,\,\Pi(\cl F.)\Omega\big)\to 0$ as $n\to\infty$,
but $\|\varphi_n\|$ does not converge to zero. (This means there is some
nonzero $\wt\varphi\in\cl H._\Gamma$ which is orthogonal to $\Pi(\cl F.)\Omega$).
Let $\varphi_n=\Pi(A_n)\Omega,$ $A_n\in\cl C.',$ then
$$\eqalignno{(\varphi_n,\,\Pi(B)\Omega)&=
\int_X\psi_x(A_n^*B)\, d\mu(x)=\big(\oo,A_n^*B\oo\big)\cr
&=\big(A_n\oo,B\oo\big)\to 0\qquad\forall\; B\in\pi_\omega(\cl F.).\cr}$$
However $\oo$ is cyclic for $\pi_\omega(\cl F.)$ in $\cl H._\omega,$ hence
$A_n\oo\to 0$ and so since
$\|\varphi_n\|^2=\|A_n\oo\|^2\to 0$, we have contradicted the
hypothesis, so $\kappa(\Omega)$ must be cyclic for $\Pi(\cl F.)$. Now since
$$(\oo,\,A\oo)=\wt\omega(A)=\int_X\psi_x(A)\, d\mu(x)= (\Omega,\,\Pi(A)\Omega)$$
the unitary equivalence follows from the GNS--theorem.

\itemitem{\bf Remarks.}{\bf (\rn1)} Since $\big(\pi\s\tilde\omega.,\,
\cl H.\s\tilde\omega.,\,\Omega\s\tilde\omega.\big)$
is the representation in which $\cl C.'$ is defined, $\Pi$ is an
isomorphism.
\itemitem{\bf (\rn2)} We can realise  vectors $\xi\in\cl H._{\Gamma}$
as  sections $x\to\xi(x)\in\cl H._{\psi_x}$ (almost everywhere)
such that 
$$(\xi,\phi)=\int_X(\xi(x),\,\phi(x))_{\cl H._{\psi_x}}d\mu(x),$$
but we will not need it here, so omit it.
\def\ap{\wt A\s\varphi,\,\zeta.}
\def\Ap{A\s\varphi,\,\zeta.}
\item{\bf Def.} Let $\varphi,\;\zeta\in\Pi(\cl C.')\Omega\subset
\Gamma(\cl H.(\psi))$ and define a section
$\ap\in\Gamma(\cl B.(\psi))$ by
$$\big(\ap c\big)(x):=\varphi(x)\,\big(\zeta(x),\, c(x)\big)\s{\cl H._{\psi_x}}.\qquad
\forall\; c\in\Gamma(\cl H.(\psi))\;.$$
\thrm Lemma 2.3."$\ap$ preserves $\Pi(\cl C.')\Omega$ and
$\ker(\cdot,\cdot)$, and is bounded on these spaces, hence
defines an operator $\Ap\in\cl B.(\cl H._\Gamma)$
by $\Ap\kappa(\eta):=\kappa(\ap\eta)$ for $\eta\in\Pi(\cl C.')
\Omega$." 
Let $\varphi=\Pi(A)\Omega$, $\zeta=\Pi(B)\Omega$, 
$\eta=\Pi(C)\Omega$, $A,\, B,\, C\in\cl C.'$, then
$$\eqalignno{\big(\ap\eta\big)(x)&=
\pi\s\psi_x.(A)\Omega\s\psi_x.\cdot\left(\pi\s\psi_x.(B)\Omega\s\psi_x.,\;
\pi\s\psi_x.(C)\Omega\s\psi_x.\right)\cr
&=\psi_x(B^*C)\cdot\pi\s\psi_x.(A)\Omega\s\psi_x.
=\Phi(\Upsilon^{-1}(\delta(B^*C)))(x)\cdot\Pi(A)\Omega(x)\cr
&=\Pi(\Upsilon^{-1}(\delta(B^*C))A)\Omega(x)\cr}$$
where we used the fact that $\Upsilon^{-1}(\delta(B^*C))\in\cl C.''=\cl C.$
and theorem 2.1. Thus $\ap\eta\in\Pi(\cl C.')\Omega$. That 
$\ap$ preserves $\ker(\cdot,\cdot)$ will follow from the next calculation as
well as the remaining claims. Since by theorem 2.2, $\kappa(\Omega)$
is cyclic for $\pi(\cl F.)$, it suffices to do the calculation
on $\Pi(\cl F.)\Omega$, so now let $\eta=\Pi(C)\Omega$, $C\in\pi_\omega(\cl F.)$.
Then
$$\eqalignno{\big\|\ap\eta\big\|^2&=\int\big|\psi_x(B^*C)\big|^2\,
\psi_x(A^*A)\, d\mu(x)\cr
&\leq\int_X\psi_x(B^*B)\cdot\psi_x(C^*C)\cdot\psi_x(A^*A)\, d\mu(x)
\qquad\qquad\cr
&\leq\|A\|^2\|B\|^2\int_X\psi_x(C^*C)\, d\mu(x)\cr
&=\|A\|^2\|B\|^2\|\eta\|^2\;.\cr}$$

\noindent Define the norm $\|a\|=\sup\limits_{x\in X}\big\|a(x)\|\s{
\cl B.(\cl H._{\psi_x})}.$ on the space of bounded sections:\chop
$\Gamma_0(\cl B.(\psi)):=
\set a\in\Gamma(\cl B.(\psi)),\|a\|<\infty.$
 which makes it into a
C*--algebra.
\thrm Lemma 2.4." $\set\ap,{\varphi,\;\zeta\in\Pi(\cl C.')\Omega}.
\subset\Gamma_0(\cl B.(\psi))$."
$(\ap c)(x)=\varphi(x)\,\big(\zeta(x),\, c(x)\big)$ for $c\in\Gamma
(\cl H.(\psi))$, hence
$$\eqalignno{\big\|(\ap\, c)(x)\big|&\leq\|\varphi(x)\|\cdot
\|\zeta(x)\|\cdot\|c(x)\|\;,\qquad\quad\hbox{and}\cr
\big\|\ap(x)\big\|\s{\cl B.(\cl H._{\psi_x})}.&=\|\varphi(x)\|\cdot
\|\zeta(x)\|\;,\qquad\qquad\hbox{so}\cr
\big\|\ap\big\|&=\sup_{x\in X}\left(\|\varphi(x)\|\cdot\|\zeta(x)\|
\right)=\sup_{x\in X}\Big(\psi_x(A^*A)\,\psi_x(B^*B)\Big)^{1/2}\cr
&\leq\|A\|\cdot\|B\|<\infty\cr}$$
where we assumed $\varphi=\Pi(A)\Omega$, $\zeta=\Pi(B)\Omega$.

\noindent Using the C*--operations of $\Gamma_0(\cl B.(\psi))$, we now define
$$\wt{\cl L.}:=C^*\set\ap,{\varphi,\,\zeta\in\Pi(\cl C.')\Omega}.
\subset\Gamma_0(\cl B.(\psi))\;.$$
\thrm Theorem 2.5." With the data $\cl F.,\;\omega,\;\cl C.=\cl C.''$ above;--\chop
$(\rn1)\quad\wt{\cl L.}(x):=\set a(x),a\in{\cl L.}.=\cl K.(\cl H.\s\psi_x.)$\chop
$(\rn2)$ $\Pi(\cl C.')\wl\subset\wl\supset\wl
\,\Pi(\cl C.')$ and $[\wl,\,\Pi(\cl C.)]=0$,\chop
$(\rn3)$ $\wl$ lifts through $\kappa$ to define a representation
\chop $\rho:\wl\to\cl B.(\cl H._\Gamma)$. Thus by $(\rn2)$ $\rho(\wl)
\subset\pi(\cl C.')$."
$(\rn1)$ In a fibre we have
$\big\|\ap(x)\big\|\s{\cl B.(\cl H._{\psi_x})}.=\|\varphi(x)\|\cdot
\|\zeta(x)\|$ so $\ap(x)$ is continuous in $\varphi(x)$, $\zeta(x)$,
and thus the norm closure of $\set\ap(x),\varphi,\;\zeta\in\Pi(
{\cl C.}')\Omega.$ contains all rank one operators, using the fact that
$\big(\Pi(\cl C.')\Omega\big)(x)$ is dense in $\cl H.\s\psi_x.$.
Since
$\cl K.(\cl H.\s\psi_x.)$ is spanned by its rank one operators
and closure in a supremum norm produces pointwise closure, we get
that 
$\wt{\cl L.}(x)=\cl K.(\cl H.\s\psi_x.)$.\chop
$(\rn2)$ Let $\varphi=\Pi(A)\Omega$, $\zeta=\Pi(B)\Omega$, then
for $E\in\cl C.'$:
$$\eqalignno{\big(\pi(E)\ap c\big)(x)&=\pi\s\psi_x.(EA)\Omega\s\psi_x.
\cdot\left(\pi\s\psi_x.(B)\Omega\s\psi_x.\,,\; c(x)\right)\cr
&=\big(\wt  A\s\xi,\zeta.\, c\big)(x)\qquad\quad\hbox{where}\;\;
\xi:=\Pi(EA)\Omega\,.\cr}$$
So $\Pi(E)\ap\in\wl$. Similarly $\ap\Pi(E)\in\wl$, and as
$\ap(x)^*=\wt A\s\zeta,\varphi.(x)$ it follows that 
$\Pi(\cl C.')\wl\subset\wl\supset\wl\,\Pi(\cl C.')$.
To see that $[\wl,\,\Pi(\cl C.)]=0$, let $F\in\cl C.$, so
$$\eqalignno{\big(\ap\Pi(F)c\big)(x)&=\varphi(x)\,\big(\zeta(x),\,
\Pi(F)c(x)\big)=\varphi(x)\,\big(\zeta(x),\,\Phi(F)(x)\,c(x)\big)\cr
&=\Phi(F)(x)\,\big(\ap c\big)(x)=\big(\Pi(F)\ap c\big)(x)\;.\cr}$$
$(\rn3)$ By lemma 2.3 we already know that $\ap$ lifts through
$\kappa$, so since lifting is a homomorphism, we define
$\rho(\ap):=\Ap\in\cl B.(\cl H._\Gamma)$ and check continuity.
Let $\varphi=\Pi(A)\Omega$, $\zeta=\Pi(B)\Omega$, $\eta=\Pi(C)\Omega$
then
$$\eqalignno{\big\|\rho(\ap)&\kappa(\eta)\big\|^2
=\big\|\kappa(\ap\eta)\big\|^2
=\int_X\big\|(\ap\eta)(x)\big\|^2d\mu(x)\cr
&=\int_X\big|\psi_x(B^*C)\big|^2\cdot\psi_x(A^*A)\, d\mu(x)\cr
&\leq\int_X\psi_x(B^*B)\cdot\psi_x(C^*C)\cdot\psi_x(A^*A)\, d\mu(x)\cr
&\leq\sup_{x\in X}\big(\psi_x(A^*A)\cdot\psi_x(B^*B)\big)\int_X
\psi_x(C^*C)\, d\mu(x)\cr
&=\big\|\ap\big\|^2\cdot\|\kappa(\eta)\|^2 \cr}$$

\itemitem{\bf Remarks.}{\bf (1)} Note that theorem 2.5(\rn2) expresses that
$\wl$ is a $C(X)\hbox{--algebra}$, given that $\Pi(\cl C.)=C(X)$.
Thus the family $\set\wl(x),x\in X.$ can be topologised as a Fell bundle (cf. 
M. Nilsen~[Ni]). In fact, since $\cl C.'$ is also obviously
a $C(X)\hbox{--algebra}$, it is also the set of continuous sections of a Fell
bundle. Note also that $\wl$ has ideals corresponding to closed subsets $Y$ of $X$
by $\cl I._Y:=\set L\in\cl L.,\Pi(L)\rest Y=0..$ Obviously, since
$\wl$ is fibrewise the compacts, it is a nonunital algebra.
\itemitem{\bf (2)} Another, possibly smaller choice for $\wl$ is
$$\cl L._{\cl F.}:=C^*\set\ap,\varphi,\;\zeta\in{\rm Span}\{\Pi({\cl C.\pi_\omega
(\cl F.)''})\Omega\}.$$
in which case we still have fibrewise the compacts,
$\cl L._{\cl F.}(x)=\cl K.(\cl H.\s\psi_x.)$, but instead of 2.4(\rn2)
we now have only the weaker property that $\Pi(\cl C.)$ and $\Pi(\pi_\omega(
\cl F.)'')$ are in the relative multiplier algebra of $\cl L._{\cl F.}$
in $\Gamma_0(\cl B.(\psi))$.\chop
\thrm Lemma 2.6."$(\rn1)$ Let $a\in\wl$, then the map $x\to\|a(x)\|$
is continuous.\chop
$(\rn2)$ Fix an $x\in X$, then $\set a\in\wl,a(x)=0.$ is dense in the set
$\set f\cdot a,a\in\wl,\; f\in C(X),\; f(x)=0.$."
$(\rn1)$ First we show that $x\to\|a(x)\|$ is continuous for the generating set
of $\wl$. Let $a=\ap$ with $\varphi=\Pi(A)\Omega$, $\zeta=\Pi(B)\Omega$,
$A,\; B\in\cl C.'$, then
$$\eqalignno{\big\|a(x)\big\|&=\left\|\ap(x)\right\|=\|\varphi(x)\|\cdot
\|\zeta(x)\|\cr
&=\left\|\pi\s\psi_x.(A)\Omega\s\psi_x.\right\|\cdot
\left\|\pi\s\psi_x.(B)\Omega\s\psi_x.\right\|
=\left[\psi_x(A^*A)\;\psi_x(B^*B)\right]^{1/2}\cr}$$
which is continuous in $x$ because 
$\psi_x(A^*A)=\Phi(\Upsilon^{-1}(\delta(A^*A)))(x)$
and this is continuous by theorem 2.1(\rn1). It is easy to see that
$x\to\|a(x)\|$ is continuous for linear combinations of the $\ap$:
let $a,\; b\in\Gamma(\cl B.(\psi))$ be such that $x\to\|a(x)\|$ and
$x\to\|b(x)\|$ is continuous, then for a convergent net $x_\nu\to x$ in $X$,
$$\eqalignno{\big|\|a(x_\nu)+b(x_\nu)\|&-\|a(x)+b(x)\|\big|
\leq\big\|a(x_\nu)+b(x_\nu)-a(x)-b(x)\big\|\cr
&\leq\|a(x_\nu)-a(x)\|+\|b(x_\nu)-b(x)\|\longrightarrow\; 0\,.\cr}$$
Next we show that $\cl Y.:={\rm Span}\set\ap,\varphi,\;\zeta\in\Pi({\cl F.})
\Omega.$ is in fact a dense *--subalgebra of $\wl$, not just a generating set.
That it is involutive follows from $\ap^*=\wt A\s\zeta,\varphi.$.
We only need to show that $\ap\cdot\wt A\s\xi,\eta.\in\cl Y.$. Let
$c\in\Gamma(\cl H.(\psi))$, then
$$\eqalignno{\big(\ap &\cdot\wt A\s\xi,\eta.\;c\big)(x)=
\varphi(x)\;\big(\zeta(x),\,\big(\wt A\s\xi,\eta.\, c\big)(x)\big)\cr
&=\varphi(x)\;\big(\zeta(x),\,\xi(x)\big)\cdot\big(\eta(x),\, c(x)\big)
=\big(\wt A\s\bar\varphi,\eta.c\big)(x)\cr}$$
where $\bar\varphi:=\Pi(\wt{(\zeta,\xi)})\varphi\in\Pi(\cl C.')\Omega$,
since $x\to\wt{(\zeta,\xi)}(x):=(\zeta(x),\,\xi(x))$ is continuous so
corresponds to an element $\wt{(\zeta,\xi)}$
 of $\cl C.$. Thus \chop$\ap\cdot\wt A\s\xi,\eta.\in\cl Y.$.
\chop Finally, we need to show that if $\{a_n\}\subset\cl Y.$ is a sequence
converging to an $a\in\wl$ in norm, then $x\to\|a(x)\|$ is continuous.
Given a converging net $x_\nu\to x$ in $X$, 
$$\eqalignno{\Big|\|a(x_\nu)\|-\|a(x)\|\Big|&\leq\|a(x_\nu)-a(x)\|\cr
&=\big\|a(x_\nu)-a_n(x_\nu)+a_n(x_\nu)-a_n(x)+a_n(x)-a(x)\big\|\cr
&\leq\|a(x_\nu)-a_n(x_\nu)\|+\|a_n(x_\nu)-a_n(x)\|+\|a_n(x)-a(x)\|\cr
&\leq 2\|a-a_n\|+\|a_n(x_\nu)-a_n(x)\|\cr}$$
and this can be made arbitrary small for suitable choices of $n$ and $\nu$.\chop
$(\rn2)$ Denote $\cl I._x:=\set a\in\wl,a(x)=0.$ and $\cl J._x:=
\chop\overline{\set f\wl,f\in C(X),\; f(x)=0.}$.
If $g=f\cdot a\in\cl J._x$, we clearly have that $g(x)=f(x)\, a(x)=0$,
i.e. $\cl J._x\subseteq\cl I._x$. Conversely, let $a\in\cl I._x$, then
by part (\rn1) of this lemma, $U_\varepsilon:=\set y\in X,\|a(y)\|<\varepsilon.$
is open for any $\varepsilon>0$. Since $X$ is the spectrum of a von Neumann 
algebra, it is completely regular, so there is a continuous function 
$f:X\to[0,\, 1]$ such that $f(x)=0$, $f(y)=1$ for all $y\not\in U_\varepsilon$.
So $$\big\|(a-fa)(x)\big\|=\|a(x)(1-f(x))\|=\|a(x)\|\cdot |1-f(x)|=0$$
as $a(x)=0$. Since $\|(a-fa)(y)\|<\varepsilon$ for all $y\in U_\varepsilon$
and \chop $\|(a-fa)(y)\|=0$ when $y\not\in U_\varepsilon$, it is clear that
$\cl I._x=\cl J._x$.

\thrm Lemma 2.7." Any pure state $\gamma$ of $\wl$ is of the form
$\gamma(a)=\gamma_x(a(x))$ for all $a\in\wl$ where $x$ is a distinguished
point $x\in X$, and $\gamma_x$ is a pure state of $\wl(x)=\cl K.(\cl H.\s\psi_x.)$."
{}From 2.5 we have $C(X)\wl\subset\wl\supset\wl\, C(X)$ using $\Pi(\cl C.)=C(X)$
in the explicit action, cf. the definition below 2.1.
Hence by Dixmier 2.11.7 [Di] there is a unique extension $\wt\gamma$ of
$\gamma$ to a pure state of \chop $C^*(\wl\cup\Pi(\cl C.))\subset M(\wl)$.
Since $C(X)=\Pi(\cl C.)\subset Z\big(C^*(\wl\cup\Pi(\cl C.)\big)$, 
$\wt\gamma\rest C(X)$ is pure (since the image of the centre
of a C*--algebra under an irreducible representation is one--dimensional).
Thus $\wt\gamma\rest C(X)$ is evaluation at some distinguished point 
$x\in X$, and obviously the left kernel
$$N\s\tilde\gamma\rest C(X).= \ker\wt\gamma\rest C(X)=\set f\in C(X), f(x)=0.\,.$$
Since $C(X)$ commutes with $\wl$, this means
$$\big|\gamma(Lf)\big|^2\leq\gamma(L^*L)\,\gamma(f^*f)=0\;\;\forall\;
f\in\ker\wt\gamma\rest C(X)\;,\;\; L\in\wl$$ i.e.
$\set\wl f, f(x)=0.\subset N_\gamma$. Since $\gamma$ is pure,
$$\ker\gamma=N_\gamma+N^*_\gamma\supset\overline{\set\wl f, f(x)=0.}$$
using Dixmier 2.9.1 [Di]. The last set is a two--sided ideal, hence in
$\ker\pi_\gamma$, and moreover by lemma 2.6(\rn2),
$\wl(x)=\wl\Big/\overline{\big\{\wl f\;\big|\; f(x)=0\big\}}$ and so
$$\gamma(a)=\gamma\big(a+\overline{\big\{\wl f\;\big|\; f(x)=0\big\}}\big)
=\hat\gamma(a(x))$$
defines a state $\hat\gamma$ on $\wl(x)$, henceforth denoted by $\gamma_x$,
i.e. $\gamma(a)=\gamma_x(a(x))$. Clearly if $\gamma_x$ is not pure, we can write it
as a convex combination of states, which through the last expression produces
a convex combination for $\gamma$, contradicting the fact that 
$\gamma$ is pure. Thus $\gamma_x$ must also be pure.

\itemitem{\bf Remarks.}{\bf (1)} Given that $\wl$ can be realised as the continuous
sections of a Fell bundle, lemma 2.7 is well--known, cf. e.g. Fell and Doran
Prop. 8.8 p582 [FD].
\def\wc{\wt{\cl C.}}
\itemitem{\bf (2)} C*--algebras of the form of $\wl$ are well--studied in the literature 
as fields of elementary algebras, cf. Dixmier [Di].
\thrm Theorem 2.8." Let $\gamma$ be a state on $\wl$, then there is a 
probability measure $\nu$ on $X$ and a $\nu\hbox{--almost}$ everywhere
defined map $\rho:\supp\nu\to\prod\limits_{x\in X}\S(\wl(x))$
such that $\rho_x\in\S(\wl(x))$, and 
$$\gamma(a)=\int_X\rho_x(a(x))\, d\nu(x)\quad\forall\; a\in\wl\;.$$"
Since by 2.5 we have $\Pi(\cl C.')$ in $M(\wl)$, both $\gamma$ and
$\pi_\gamma$ extend uniquely (on the same space) to it.
Consider the unital C*--algebra $\wc:=\pi_\gamma(\Pi(\cl C.))\subset
\pi_\gamma(\wl)'$. It will only be a von Neumann algebra if $\gamma$
is normal on $\Pi(\cl C.)$, which we cannot assume. Nevertheless,
we can apply theorem 2.1 (which works for also for a commutative C*--algebra
in the commutant) to the triple $\cl L.,\;\gamma,\;\wc$
to obtain a probability measure $\wt\nu$ on the spectrum $Y$ of
$\wc$ with support $Y$, and the Gel'fand isomorphism
$\wt\Phi:\wc\to L^\infty(Y,\,\wt\nu)$ and a $\wt\nu\hbox{--a.e.}$
defined map $\wt\psi:Y\to\S(\wl)$ such that
$$\gamma(a)=\int_Y\wt\psi_x(a)\, d\wt\nu(x)\qquad\hbox{and}\qquad
\wt\psi_x(C.a)=\wt\Phi(C)(x)\cdot\wt\psi_x(a)$$
for all $a\in\wl$, $C\in\wc$.
Now $\wc$ is a homomorphic image of $\cl C.\cong C(X)$, and as all ideals in 
$C(X)$ are of the form $\set f\in C(X),f(Z)=0.$ for some closed set $Z\subset X$,
the homomorphic images of $C(X)$ are all isomorphic to the algebras
$C(\overline{X\backslash Z})=C(V)$ for $V\subset X$ the closure of an open set.
Thus there is a homeomorphism $\beta:Y\to V\subset X$, with
$V$ the closure of an open set
and we obtain it as follows.
\def\wi{\widetilde{\cl I.}}
Since $\pi_\gamma$ must map maximal ideals of $\cl C.$ to maximal ideals
of $\wc$, define $\beta:Y\to X$ by $\pi_\gamma(\cl I.\s\beta(x).)=\wi_x$
where $\cl I._x:=\set C\in\cl C.,\Phi(C)(x)=0.$ and
$\wi_x:=\set D\in\wc,\wt\Phi(D)(x)=0.$.
Observe that this implies
$$\wt\Phi(\pi_\gamma(C))(x)=\Phi(C)(\beta(x))\;.$$
We can thus write the integral of $\gamma$ over $X$ instead of $Y$;-
define a measure $\nu$ on $X$ by setting it equal to $\nu
=\wt\nu\circ\beta^{-1}$ on $\beta(Y)\subset X$, and zero outside of
this set. Then
$\gamma(a)=\int_X\wt\psi_{\beta^{-1}(x)}(a)\, d\nu(x)$ for $a\in\wl$.
In order to define the map $\rho_x$ from this, we need to show for
$\nu\hbox{--almost}$ all $x$ that $\wt\psi_{\beta^{-1}(x)}(a)$
only depends on the value $a(x)$ of $a\in\wl$. Now
$$\wt\psi_{\beta^{-1}(x)}(C.a)=\wt\Phi(C)(\beta^{-1}(x))\cdot
\wt\psi_{\beta^{-1}(x)}(a)$$
so if $C\in\wi_{\beta^{-1}(x)}=\pi_\gamma(\cl I._x)$, then 
$\wt\psi_{\beta^{-1}(x)}(C.a)=0$ for all $a\in\wl$.
But by lemma 2.6(\rn2),
$$\overline{\cl I._x\cdot\wl}=\set a\in\wl,a(x)=0.=:K_x\;$$
so $\wt\psi_{\beta^{-1}(x)}(K_x)=0$, and thus using
$\wl(x)=\wl/K_x$, we have
$$\wt\psi_{\beta^{-1}(x)}(a)
=\wt\psi_{\beta^{-1}(x)}(a+K_x)=:\rho_x(a(x))\;,\quad a\in\wl$$
and obviously $\rho_x\in\S(\wl(x))$.
Thus finally, $\gamma(a)=\int_X\rho_x(a(x))\; d\nu(x)$ for $a\in\wl$.

\thrm Corollary 2.9." With the data $\cl F.,\;\omega,\;\cl C.$ above,\chop
$(\rn1)$ for any state $\gamma$ of $\wl$ there is a probability measure
$\nu$ on $X$ and a $\nu\hbox{--almost}$ everywhere defined map
$T:\supp\nu\to\prod\limits_{x\in X}\cl K.(\cl H.\s\psi_x.)$
such that $T_x\in\cl K.(\cl H.\s\psi_x.)$ is positive, trace class,
normalised, and satisfies
$$\gamma(a)=\int_X{\rm Tr}\big(a(x)\, T_x\big)\; d\nu(x)\;,\qquad\forall\;
a\in\wl\,.$$
$(\rn2)$ Let $\theta(\gamma)$ denote the unique extension of a state
$\gamma\in\S(\wl)$ to $\cl C.'$, then
$$\theta(\gamma)(A)=\int_X{\rm Tr}\big(\pi\s\psi_x.(A)\, T_x\big)\; d\nu(x)
\;,\qquad\forall\; A\in\cl C.'$$
where $\nu$ and $T$ are as above."
$(\rn1)$ This follows directly from theorem 2.8 and the fact that
$\wl(x)=\cl K.(\cl H.\s\psi_x.)$ for all $x\in X$.
So, since all states on the compacts are of the form
$\phi(A)={\rm Tr}(AT)$ with $T$ positive trace--class, it follows that
$\rho_x(a(x))={\rm Tr}\big(a(x)\, T_x\big)$.\chop
$(\rn2)$ Let $\gamma\in\S(\wl)$ be as above in $(\rn1)$, i.e.
 $$\gamma(a):=\int_X{\rm Tr}\big(a(x)\, T_x\big)\; d\nu(x)$$
where $a\in\wl$. Now we know the unique extension
of $\gamma$ to $\cl C.'$ is given by
$$\wt\gamma(a)=\lim_\alpha\gamma(E_\alpha a)\;,\quad
a\in\cl C.'.$$
In particular, for the states $\gamma_x\in\S(\wl)$ given by
$\gamma_x(a):={\rm Tr}\big(a(x)T_x),$ we get
$$\wt\gamma_x(a)=\lim_\alpha
{\rm Tr}\big(E_\alpha(x)\, a(x)\, T_x\big)
={\rm Tr}\big(a(x)T_x)$$
since $\gamma_x$ is normal, so has a unique extension by the last
formula to $\cl B.(\cl H._{\gamma_x}).$
Furthermore, $\|\wt\gamma_x\|=1$ for all $x,$ and as $\nu$ is a probability
measure on $X,$ the function $1$ is in $L^1(X,\,\nu).$
Thus we may apply the Lebesgue dominated convergence theorem:
$$\leqalignno{\wt\gamma(a)&=\lim_\alpha\int_X
{\rm Tr}\big(E_\alpha(x)\, a(x)\, T_x\big)\, d\nu(x) \cr
&=\int_X\lim_\alpha{\rm Tr}\big(E_\alpha(x)\, a(x)\, T_x\big)\, d\nu(x) \cr
&=\int_X{\rm Tr}\big(a(x)\,T_x)\, d\nu(x)\,. \cr}$$

\thrm Theorem 2.10." $\gamma$ is a normal state of $\cl C.'$ or of
$\pi_\omega(\cl F.)''\subset\cl C.'$ iff it can be written
$$\gamma(A)=\int_X{\rm Tr}\,\big(\pi\s\psi_x.(A)\, T_x\big)\cdot
f(x)\, d\mu(x)\;,\quad A\in\cl C.'$$
where $T_x\in\cl K.(\cl H.\s\psi_x.)$ is a.e. trace--class, positive
and normalised, and $f\in L_+^1(X,\,\mu)$ with $\mu$ the measure 
associated to the initial choice $\omega$ and $\cl C.$."
Let $\gamma$ be a normal state on $\cl C.'$ or $\pi_\omega(\cl F.)''$,
then by theorem 2.2 there is a unitary $U:\cl H._\omega\to\cl H._\Gamma$
intertwining $\pi_\omega$ with $\Pi$. Thus by Kadison and Ringrose 
7.1.12 [KR], there is a countable set of vectors
$\{\varphi_n\}\subset\cl H._\Gamma$ such that $1=\sum\limits_{n=1}^\infty
\|\varphi\|^2$ and $\gamma(A)=\sum\limits_n(\varphi_n,\,\Pi(A)\varphi_n)$.
Now let $\{\zeta_k^n\}\subset\Pi(\cl C.')\Omega$ be sequences such that
$\kappa(\zeta_k^n)\to\varphi_n\in\cl H._\Gamma$ where the convergence is in
$k$. We can in fact choose such sequences with $\|\kappa(\zeta_k^n)\|
=\|\varphi_n\|$ because if $\zeta_k^n$ is a nonzero sequence converging
to $\varphi_n$, so is $\zeta_k^n\|\varphi_n\|\big/\|\zeta_k^n\|$.
Below we will blur the distinction between $\kappa(\zeta_k^n)$ and
$\zeta_k^n$. Thus
$$\eqalignno{1=\sum_{n=1}^\infty\|\varphi_n\|^2&=
\sum_n\|\zeta_k^n\|^2=\sum_n\int_X\big\|\zeta_k^n(x)\big\|^2\s
{\cl H.}_{\psi_x}.\; d\mu(x)\cr
&=\int_X\sum_n
\big\|\zeta_k^n(x)\big\|^2\s
{\cl H.}_{\psi_x}.\; d\mu(x)\cr}$$
by Fubini and absolute convergence. Thus
$$\leqalignno{\sum_n
\big\|\zeta_k^n(x)\big\|^2\s
{\cl H.}_{\psi_x}.&\in L^1(X,\,\mu)_+\;.\qquad\hbox{Now}\cr
\Big|\big(\zeta_k^n(x),\, B\,\zeta_k^n(x)\big)\s{\cl H.}_{\psi_x}.\Big|&
\leq\|B\|\cdot\big\|\zeta_k^n(x)\big\|^2\s{\cl H.}_{\psi_x}.
\qquad\hbox{hence}\cr
\gamma_x^k(B):=\sum_n&
\big(\zeta_k^n(x),\, B\,\zeta_k^n(x)\big)\s{\cl H.}_{\psi_x}.\cr}$$
defines a positive functional on $\cl B.(\cl H.\s\psi_x.)$ for 
$\mu\hbox{--almost}$ all $x$. Moreover, it is normal by Kadison and Ringrose
7.1.12 [KR]. Now
for $A$ a positive element of $\cl C.'$
$$\eqalignno{\gamma(A)&=\sum_n(\varphi_n,\, A\varphi_n)=\sum_n\lim_k
(\zeta_k^n,\, A\zeta_k^n)=
\lim_k\sum_n
(\zeta_k^n,\, A\zeta_k^n)\cr}$$
by dominated convergence, since $\big|(\zeta_k^n,\, A\zeta_k^n)\big|
\leq\|A\|\cdot\|\zeta_k^n\|^2=\chop
\|A\|\cdot\|\varphi_n\|^2$. So
$$\eqalignno{\gamma(A)
&=\lim_k\sum_n\int_X\left(\zeta_k^n(x),\; A(x)\,\zeta_k^n(x)\right)\s{
\cl H.}_{\psi_x}.\; d\mu(x)\cr
&=\lim_k\int_X\sum_n\left(\zeta_k^n(x),\; A(x)\,\zeta_k^n(x)\right)\s{
\cl H.}_{\psi_x}.\; d\mu(x)\cr
&=\lim_k\int_X\gamma_x^k(A(x))\; d\mu(x)\cr}$$
where we used Fubini's theorem and absolute convergence.
Moreover, since
$$\big|\gamma_x^k(A(x))\big|\leq\|A\|\;\sum_n\big\|\zeta_k^n(x)\big\|^2\s{
\cl H.}_{\psi_x}.\in L^1(X,\,\mu)_+$$
we can use the dominated convergence theorem to conclude
$$\gamma(A)=\int_X\lim_k\gamma_x^k(A(x))\; d\mu(x)\eqno{(1)}$$
providing we can show that the pointwise limits
$\lim\limits_k\gamma_x^k(A(x))$ exist a.e. which is what we prove now.
Since $\zeta_k^n\in\Pi(\cl C.')\Omega$, let $\zeta_k^n=\Pi(A_k^n)\Omega$,
$A_k^n\in\cl C.'$, so for $B\in\cl C.'$
$$\eqalignno{\gamma_x^k(B)&=\sum_n\big(\zeta_k^n(x),\, B(x)\,\zeta_k^n(x)\big)
=\sum_n\big(\Pi(A_k^n)\Omega(x),\,\Pi(BA^n_k)\Omega(x)\big)\cr
&=\sum_n\psi_x({A_k^n}^*B\, A_k^n)\;\qquad\quad\hbox{so}\cr
\big|\gamma_x^k(B)-&\gamma_x^\ell(B)\big|=\Big|\sum_n\psi_x({A_k^n}^*B\, A^n_k
-{A^n_\ell}^*B\, A_\ell^n)\Big|\cr
&\leq\Big|\sum_n\big(
\psi_x((A_k^n-A_\ell^n)^*B(A_k^n-A_\ell^n))
+\psi_x({A_\ell^n}^*B(A^n_k
-{A^n_\ell}))\cr &\qquad\qquad+\psi_x((A_k^n-A^n_\ell)^*BA^n_\ell)\big)\Big|\cr
&\leq\|B\|\sum_n\big(
\psi_x((A_k^n-A_\ell^n)^*(A_k^n-A_\ell^n))
\cr &\quad+ 2\psi_x({A_\ell^n}^*A^n_\ell)^{1/2}
\cdot\psi_x((A_k^n-A^n_\ell)^*(A_k^n-A^n_\ell))^{1/2}\big)&(2)\cr}$$
using the Cauchy--Schwartz inequality. Now
$$\sum_n
\psi_x((A_k^n-A_\ell^n)^*(A_k^n-A_\ell^n))
=\left\|\zeta_k^n(x)-\zeta_\ell^n(x)\right\|^2$$
which must converge to zero a.e. as $k$ and $\ell$ approach infinity because
$$\left\|\zeta_k^n-\zeta_\ell^n\right\|=\Big(\int_X
\left\|\zeta_k^n(x)-\zeta_\ell^n(x)\right\|^2d\mu(x)\Big)^{1/2}$$
and this approaches zero with $k$, $\ell$ since $\{\zeta_k^n\}$
is a convergent sequence. Applying the Cauchy--Schwartz inequality to
the sum:
$$\eqalignno{
\Big|\sum_n
\psi_x({A_\ell^n}^*A^n_\ell)^{1/2}
&\cdot\psi_x((A_k^n-A^n_\ell)^*(A_k^n-A^n_\ell))^{1/2}\Big|^2
\cr &\leq\sum_n\psi_x({A_\ell^n}^*A_\ell^n)\cdot
\sum_m\psi_x((A_k^m-A_\ell^m)^*(A_k^n-A_\ell^n))\cr}$$
which as we saw approach zero a.e. Thus (2) converges to zero
a.e. and so the pointwise limit $\lim_k\gamma_x^k(B(x))$
exists a.e. and (1) is justified.\chop
Since the normal states are sequentially weak*--closed
(cf. Takesaki 5.2 p148 [Tak]), we conclude that $\lim\limits_k\gamma_x^k$
is a normal state on $\cl B.(\cl H.\s\psi_x.)$.
Thus there is a positive normalised trace--class operator $T_x
\in\cl K.(\cl H.\s\psi_x.)$ such that
$\gamma_x(A(x))=f(x)\cdot{\rm Tr}\big(A(x)\, T_x\big)$
where $f(x)=\gamma_x(\un)>0$. That $f\in L^1(X,\,\mu)$ follows from
$$\displaylines{1=\gamma(\un)=\int_X\gamma_x(\un)\, d\mu(x)=
\int_X f(x)\, d\mu(x)\;.\cr
\hbox{Thus}\qquad\qquad\gamma(A)=\int_X{\rm Tr}\big(\pi\s\psi_x.(A)\, T_x\big)
\; f(x)\, d\mu(x)\qquad\forall\; A\in\cl C.'\;.\cr}$$
Conversely, assume $\gamma(A)$ to have this form. Then by Kadison and
Ringrose 7.1.12 [KR], it suffices to show $\gamma$ is strong operator continuous
on the unit ball of $\cl C.'$. Let $\{A_\alpha\}\subset\cl C.'$ be a net
converging to zero in strong operator topology
and with $\|A_\alpha\|\leq 1$. That is, for all $\varphi\in
\cl H._\omega$
$$\|A_\alpha\varphi\|^2=\int_X\big\|A_\alpha(x)\,\varphi(x)\big\|^2\s
{\cl H.}_{\psi_x}.\; d\mu(x)\longrightarrow 0$$
which implies that $\|A_\alpha(x)\,\varphi(x)\|\to 0$ almost everywhere.
Since $A_\alpha(x)=\pi\s\psi_x.(A_\alpha)$,
$$\eqalignno{\big|\gamma(A_\alpha)\big|=\big|\int_X{\rm Tr}\left(\pi\s\psi_x.(A_\alpha)
\, T_x\right)\; f(x)\, d\mu(x)\big|\leq
\int_X\big|{\rm Tr}\left(\pi\s\psi_x.(A_\alpha)
\, T_x\right)\big|\; f(x)\, d\mu(x)\,.\cr}$$
Moreover $\big|{\rm Tr}\left(\pi\s\psi_x.(A_\alpha)
\, T_x\right)\big|\; f(x)\leq\|A_\alpha\|\cdot f(x)$
which is of course an $L^1\hbox{--function}$ and 
$\lim\limits_\alpha\big|{\rm Tr}\left(\pi\s\psi_x.(A_\alpha)
\, T_x\right)\big|\to 0$ because 
$\big\|\pi\s\psi_x.(A_\alpha)\big\|\leq 1$
and ${\rm Tr}(\cdot\, T_x)$ is a normal functional on
$\cl B.(\cl H.\s\psi_x.)$ (recall that $A_\alpha(x)$ converges
almost everywhere to $0$ in the strong operator topology).
Thus by the dominated convergence theorem,
$$\eqalignno{\lim_\alpha\big|\gamma(A_\alpha)\big|&\leq\lim_\alpha
\int_X\big|{\rm Tr}\big(\pi\s\psi_x.(A_\alpha)\, T_x\big)\big|\;
f(x)\, d\mu(x)\cr
&=\int_X\lim_\alpha\big|{\rm Tr}\big(\pi\s\psi_x.(A_\alpha)\, T_x\big)\big|\;
f(x)\, d\mu(x)=0\;.\cr}$$
The argument generalises to nets $A_\alpha\to A\not=0$ in strong operator
 topology for $\|A_\alpha\|\leq 1$, hence $\gamma$ is normal.

\itemitem{\bf Remark.} Observe that the theorem automatically
constructs a state on $\cl C.'$ even if one starts from a state
on $\pi_\omega(\cl F.)''$, though due to the choices involved
this need not be unique.
\item{\bf Def.} Given $\cl F.,\;\omega,\;\cl C.$ as above, define the set
$\S_\omega\subset\S(\wl)$ as the set of those states $\gamma$ with 
$\gamma(a)=\int_X{\rm Tr}\big(a(x)\, T_x\big)\, d\nu(x)$
as in 2.9,
where $\nu$ is absolutely continuous with respect to $\mu$.
As before, we denote the set of normal states of a concrete
von Neumann algebra $\cl N.$ by $\S_N(\cl N.)$.
\thrm Theorem 2.11." Assume the data and notation of Corr. 2.9, then\chop
$(\rn1)$ For a state $\gamma\in\S(\wl)$ the extension $\theta(\gamma)$
is normal on $\cl C.'$ or on $\pi_\omega(\cl F.)''$ iff its measure $\nu$
on $X$ is absolutely continuous with the measure $\mu$ associated with
$\omega$. \chop
$(\rn2)$ For each normal state $\eta$ of $\cl C.'$ or $\pi_\omega(\cl F.)''$
there is a $\gamma\in\S(\wl)$ such that $\eta=\theta(\gamma)$.
\chop
$(\rn3)$ We have:
 $\theta:\S_\omega\to\S_N(\pi_\omega(\cl F.)'')$ is a surjection.
  Since $\wl\subset\cl C.'$, we see that
$\theta:\S_\omega\to\S_N(\cl C.')$ is a bijection."
$(\rn1)$ By Corr. 2.9(\rn2) and theorem 2.10, this follows immediately.
\chop $(\rn2)$
Now $\cl C.'\subset\wl'',$ so normal states of $\cl C.'$
are restrictions of normal states of $\wl'',$ and these
in turn are the unique extensions of states from $\wl$
by strong operator continuity. 
Since $\cl C.'\subset M(\wl)\subset\wl'',$
these strong operator extensions are just the unique extensions by $\theta.$
Thus we have surjectivity as claimed.
\chop $(\rn3)$ This is just a restatement of the preceding parts.

\itemitem{\bf Remark.} Since $\S_\omega$ is a proper subset of
$\S(\wl),$ this means that $\wl$ cannot in general be 
an ideal host for $\cl C.'$ or for $\pi_\omega(\cl F.)'',$ a fact which is also
obvious from the commutative situation:
$\pi_\omega(\cl F.)''=L^\infty(X,\,\mu)=\wl.$
But we have almost finished showing that it is a quasi--host.
\thrm Theorem 2.12."
Let $\cl C.\cong L^\infty(X,\mu)$ be maximally commutative in $\pi_\omega(\f)'.$
Then $\cl C.'=(\cl N.\cup\cl C.)''\supset\wl,$
$\pi_\omega(\f)\subset M(\wl),$
$$\S_\omega=\set\varphi\in\S(\wl),\wt\varphi\restriction L^\infty(X,\mu)\;
\hbox{is normal}.,$$
${\wt\S_\omega\restriction\pi_\omega(\f)}=\S_N(\cl N.)=
\S_N(\pi_\omega(\f))$ and 
${\wt\S_\omega\restriction C^*(\cl N.\cup\cl C.)}$ is in bijection with
$\S_\omega.$ In short, $\wl$ is a quasi--host for
$(\cl N.,\S_N(\cl N.))$ (hence for $(\pi_\omega(\f),\S_N(\pi_\omega(\f))).$"
Since $\cl C.$ is maximally commutative in $\cl N.',$ we have
$\cl C.=\cl N.'\cap\cl C.'=(\cl N.\cup\cl C.)'$ and thus
$\cl C.'=(\cl N.\cup\cl C.)''$ hence a normal state on $\cl C.'$
is uniquely determined by its values on $\cl C.$ and on $\cl N..$
Thus by Theorem~2.11(\rn3) we have that
$\wt\S_\omega\restriction C^*(\cl N.\cup\cl C.)$ is in bijection with
$\S_\omega.$ 
We already have the embeddings stated, so to check the claimed characterisation of 
$\S_\omega,$ recall that it consists of states $\gamma\in\S(\wl)$ such that
$$\gamma(A)=\int_X{\rm Tr}\left(A(x)T_x\right)\,  d\nu(x)\,,\quad
A\in\wl$$
where $\nu$ is absolutely continuous w.r.t. $\mu.$
Then for $f\in\cl C.=L^\infty(X,\mu)$ we have for any approximate identity
$\{E_\alpha\}$ of $\wl:$
$$\wt\gamma(f)=\lim_\alpha\gamma(f\cdot E_\alpha)=\lim_\alpha
\int_X f(x)\,{\rm Tr}\left(E_\alpha(x)T_x\right)\,  d\nu(x)
=\int_Xf(x)\,d\nu(x)$$
using the argument in the proof of Corollary~2.9(\rn2).
These are precisely the normal states of 
$L^\infty(X,\mu).$ This completes the proof.

Return now to the original pair $(\f,\S_0)$ at the start of the investigation,
which was equivalent to the examination of $(\pi_{\S_0}(\f),\,\S_N(\pi_{\S_0}(\f))),$
where $\pi_{\S_0}=\bigoplus\limits_{\omega\in\S_0}\pi_\omega.$
Take the direct sum $\wl_0:=\bigoplus\limits_{\omega\in\S_0}\wl_\omega
\subset\cl B.(\cl H._{\S_0})$ where $\wl_\omega$ is the quasi--host
constructed above for the pair
 $(\pi_\omega(\f),\S_N(\pi_\omega(\f))).$
Then $$\pi_{\S_0}(\f)\subset M(\wl_0),\quad\hbox{and}\quad\bigoplus\limits_{\omega\in\S_0}
L^\infty(X_\omega,\mu_\omega)\subset ZM(\wl_0).$$
Let $(X,\mu)$ be the disjoint union of the measure spaces $(X_\omega,\mu_\omega)$
(hence $\mu$ is not a probability measure), so we can write
$L^\infty(X,\mu)=\bigoplus\limits_{\omega\in\S_0}L^\infty(X_\omega,\mu_\omega),$ and 
thus $L^\infty(X,\mu)\subset ZM(\wl_0).$ Moreover with notation
$$\S_\mu:=\set\varphi\in\S(\wl_0),\wt\varphi\restriction L^\infty(X,\mu)\;
\hbox{is normal}.$$
we see that for $\varphi\in\S_\mu$ that
 $\wt\varphi\restriction L^\infty(X_\omega,\mu_\omega)$ is normal
for all $\omega\in\S_0,$ hence $\S_\mu$ is the norm closed convex hull of
$\bigcup\limits_{\omega\in\S_0}\S_\omega$ and by Theorem~2.12
${\wt\S_\mu\restriction\pi_{\S_0}(\f)}=\S_0$ and
${\wt\S_\mu\restriction C^*(\pi_{\S_0}(\f)\cup L^\infty(X,\mu))}$ is in bijection
with $\S_\mu.$ In other words, $\wl_0$ is a quasi--host for $(\f,\,\S_0).$
Note that whilst $\wl_0$ is in $\cl B.(\cl H._{\S_0}),$ it need not be
in $\pi_{\S_0}(\f)'',$ (unlike ideal hosts) because
$\cl C.$ may have a part outside 
 $\pi_{\S_0}(\f)''.$

\thrm Theorem 2.13." Given the preceding notation, if $(\cl F.,\,\S_0)$
has an ideal host, then $\mu$ must have discrete points."
If $(\cl F.,\,\S_0)$ has an ideal host, then by Corollary~1.7,
$\S_0$ has pure states. So in the direct sum
$\pi_{\S_0}(\cl F.)=\bigoplus\limits_{\omega\in\S_0}\pi_\omega(\cl F.)$
we have that for some $\omega\in\S_0,$ that
$\pi_\omega(\cl F.)'=\C\un,$ and so $\cl C._\omega=\C\un$ and thus
${(X_\omega,\,\mu_\omega)}$ is the discrete trivial measure space
${(\{x\},\,\delta)},$ $\delta(\{x\})=1,$ $\delta(\emptyset)=0.$
Since this is a summand of $L^\infty(X,\mu),$
we conclude that $\mu$ has discrete points.

\noindent
This makes now explicit the sense in which we meant that $L^\infty(X,\mu)$
with $\mu$ continuous,
is an obstruction to the existence of an ideal host.
It would be nice to get a converse, i.e. to argue that an ideal
host exists iff $\mu$ has some specific structure, but we do not
have this yet.

\beginsection 3. Applications.

Here we will do a couple of applications of ideal hosts, mainly following
the well--trodden path of group algebras.
First, there is the question of useful decompositions of states
in $\S_0$ into other states in $\S_0:$
\thrm Theorem 3.1."Let $\l$ be an ideal host for a pair ${(\f,\,\S_0),}$
and let $\varphi\in\S_0.$ Then there is an integral decomposition
$$\varphi(F)=\int_{\S_0}\omega(F)\,d\mu(\omega),\quad F\in\f,$$
such that $\mu$ is pseudosupported on the pure states of
$\S_0.$ If $\l$ is separable, $\mu$ is actually supported on
the pure states of $\S_0.$"
Most of the work for the proof has already been done in the
last section.
First observe that from Theorem~1.5 we know that 
$\l''=\pi_{\S_0}(\f)'',$ and hence for an $\omega\in\S(\l)$
we have that $\pi_\omega(\l)''=\pi_{\theta(\omega)}(\f)''.$
For a given $\varphi\in\S_0,$ choose $\omega=\theta^{-1}(\varphi).$
From Theorem~2.1, for a choice of commutative algebra
$\cl C.\subset\pi_\omega(\l)',$ we have the decomposition
$$\wt\omega(A):=(\Omega_\omega,\, A\Omega_\omega)
=\int_X\psi_x(A)\, d\mu(x)$$
for all $A\in\cl C.'\supset\pi_{\wt\omega}(\f)=\pi_\varphi(\f),$
and hence by restriction to $\pi_\varphi(\f)$ this is also
a decomposition for $\varphi=\wt\omega\restriction\f\in\S_0.$
Make the usual identification of $X$ with a subset
of $\S(\l)$ by $x\to\psi_x,$ and maintain the same symbol
for the measure $\mu$ carried by this identification to $\S(\l).$
Of course via $\theta$ we can now carry this measure to $\S_0$ too.
Thus we have
$$\leqalignno{\omega(L)&=
\int_{\S(\l)}\gamma(L)\, d\mu(\gamma)\;,\qquad L\in\l,\quad\hbox{hence:}\cr
\varphi(F)&=\theta(\omega)(F)=\lim_\alpha
\int_{\S(\l)}\gamma(F\cdot E_\alpha)\, d\mu(\gamma)\cr
&=\int_{\S(\l)}\theta(\gamma)(F)\, d\mu(\gamma)
=\int_{\theta(\S(\l))}\psi(F)\, d\mu(\psi)\cr
&=\int_{\S_0}\psi(F)\, d\mu(\psi)\cr}$$
for all $F\in\f,$ any approximate identity $\{E_\alpha\}$
of $\l,$
and where we used the argument in the proof of
Corollary~2.9(\rn2) to bring the limit into the integral.
This then provides the basic decomposition formula. 
Now, since the measure $\mu$ on $\S(\l)$ is the same one
constructed in the decomposition theory in Takesaki~[Tak],
we can use his Theorem~6.28, p246, to conclude that
if $\cl C.$ is maximally commutative in $\pi_\omega(\l)',$
then $\mu$ is pseudosupported by the pure states of $\l,$
and supported by them if $\l$ is separable.
Because $\theta$ restricts to a bijection between the
pure states on $\l$ and the pure states of $\S_0,$
this proves the assertion.

\item{\bf Remark.} Since by Theorem~1.5 we deduce that $\theta$
also restricts to a bijection between the factor states of $\l$
and the factor states of $\S_0,$ we can use the central
decomposition of a state on $\l$ to obtain a similar 
decomposition of a state in $\S_0,$ in terms of factor states.
Again, we will have the two claims for the measure;--
in general it is pseudosupported on the factor states, but
if $\l$ is separable, it is actually supported by the factor states.

\medskip
As a second application for host algebras, we mention that of
inducing representations. For group algebras, this is of course
one of the main applications (cf. Rieffel~[Ri]). What allows one
to do the same here, is the relationship given in Theorem~1.9
between the representations of $\f$ and those of a host algebra.
To be more specific, given two pairs $(\f_i,\S_{0i}),$ $i=1,\; 2$
together with respective hosts $\l_i$ then providing one has
constructed a
right\ab$\l_1\hbox{--rigged}$left\ab$\l_2\hbox{--module}$
$\cl M.,$ then we can induce a representation $\pi\in{\rm Rep}\s\S_{01}.\f_1$
to a representation
$\rho\in{\rm Rep}\s\S_{02}.\f_2$ by using the map $\Lambda$ in
Theorem~1.9 to identify $\pi$ with a representation of $\l_1,$
inducing this representation via $\cl M.$ to $\l_2,$
and then identifying the result with a representation
$\rho\in{\rm Rep}\s\S_{02}.\f_2$ with  the map $\Lambda.$
The benefit of doing an induction via host algebras
 is that we remain 
within the class of representations normal w.r.t.
the representations $\pi\s\S_{0i}.$ whereas induction between
the $\f_i\hbox{'s}$ directly, can move us out of these classes.
The crux comes of course with the construction of the module
$\cl M..$ For concrete examples of this, we refer to any
example of induction via group algebras (cf. Rieffel~[Ri]).

\beginsection Acknowledgements.

I am deeply grateful to Detlev Buchholz, who has been as much responsible for the
development of host algebras as I have. Not only did the current line of thought develop out 
of his question at a seminar I gave in 1997
and many subsequent lively discussions, but he also read
several previous misguided attempts, discovering
deeply buried and serious errors.
Moreover, over the years he prodded me to return to
this project, and most recently at G\"ottingen 
provided the opportunity and support for me to pursue it, of which this
paper is the result. Thank you Detlev!

I am also grateful to the Erwin Schr\"odinger Institute in Vienna
(and in particular to Jakob Yngvason and Heide Narnhofer)
who provided me with the opportunity of completing the last part
of the paper.

\beginsection Bibliography.

\item{[AP]} C. Akemann, G.K. Pedersen: Complications of
semicontinuity in C*--algebra theory, Duke Math. J.  {\bf 40}, 785--795 (1973)
\item{[BR]} O. Bratteli, D. Robinson,
Operator Algebras and Quantum Statistical Mechanics I. Springer, New York,
1979.
\item{[Di]} J. Dixmier; C*--algebras, North--Holland,Amsterdam, 1977.
\item{[FD]} J.M.G. Fell, R.S. Doran, Representations of *-algebras, 
locally compact groups, and Banach *- algebraic bundles. Vol. 1~\&~2,
Academic Press, 1991.
\item{[GH]} H. Grundling, C.A. Hurst;
Algebraic quantization of systems with a
  gauge degeneracy. Commun. Math. Phys. {\bf 98}, 369--390 (1985)\chop
Grundling, H.: Systems with outer constraints. Gupta--Bleuler
  electromagnetism as an algebraic field theory. Commun. Math. Phys.
 {\bf 114}, 69--91 (1988)
\item{[Gr]} H. Grundling, A Group algebra for inductive limit groups.
Continuity Problems of the Canonical Commutation Relations.
Acta Applicandae Math. {\bf 46}, 107--145, 1997. 
\item{[Ha]}
Haag, R.: Local Quantum Physics. Berlin: Springer Verlag 1992
\item{[HKK]} R. Haag, R.V. Kadison, D. Kastler,
Nets of C*--algebras and Classification of States.
Commun. Math. Phys. {\bf 16}, 81--104 (1970)
\item{[KR]}
Kadison, R.V., Ringrose, J.R.: Fundamentals of the Theory of Operator Algebras II.
  Orlando: Academic Press 1986
\item{[Ku]}
Kusuda, M.: Unique state extension and hereditary C*--subalgebras.
Math. Ann. {\bf 288}, 201--209 (1990)
\item{[Ma]}
Manuceau, J.: C*--alg\`ebre de relations de commutation. Ann.
  Inst. H. Poincar\'e {\bf 8}, 139--161 (1968)
\item{[Mu]}
Murphy, G.J.: C*--Algebras and Operator Theory. Boston:
  Academic Press 1990
\item{[Ni]} Nilsen, M.: C*--bundles and $C_0(X)\hbox{--algebras.}$
Indiana Univ. Math. J., {\bf 45}, 463--477 (1996)
\item{[Pe]} 
Pedersen, G.K.: C*--Algebras and their Automorphism Groups.
  London: Academic Press 1989
\item{[Ri]} M.A. Rieffel: Induced representations of C*--algebras.
Adv. Math. {\bf 13}, 176--257 (1974)
\item{[Se]} I.E. Segal: Representations of the canonical
commutation relations, Carg\'ese lectures in theoretical physics,
Gordon and Breach, 1967.
\item{[Ta]} A. Taylor: Introduction to Functional Analysis.
John Wiley \& Sons, New York 1958.
\item{[Tak]}
Takesaki, M.: Theory of operator algebras I. New York: Springer-Verlag
  1979

\bye
\medskip
In physical situations, one usually has some distinguished group action,
e.g. for symmetries or dynamics. So, a very natural question to ask,
is whether such a group action can lift to an ideal host. To be precise,
if we have an ideal host $\l$ for a pair $(\cl F.,\,\S_0)$
and a group action $\alpha:G\to\aut\f$ such that the transposed 
action $\alpha':G\to A(\S(\f))$ preserves $\S_0,$ then does
this determine an automorphism on $\l$ compatible with the original action?
By compatible we mean the unique extension of this action on $\l$ to
$M(\l)$ agrees with $\alpha$ on $\f.$ The symbol $A(\S(\f))$
here denotes the affine w*-continuous homeomorphisms of $\S(\f).$
Unfortunately, we can only offer partial answers to this question.

First, we show how to obtain from a given automorphism $\alpha\in\aut\f$
such that $\alpha'(\S_0)=\S_0,$ a suitable automorphism of $\l''.$
That the double transpose $\alpha''$ defines an automorphism on $\f''$
is standard (cf. 7.4.5~[Pe]). Furthermore
$$\alpha'(\S_0)=\S_0=\set\omega\in\S(\f),\wt\omega(P)=1.$$
where $\wt\omega$ denotes the unique normal extension of $\omega$
to $\f''.$ Thus\chop
\centerline{$\wt{\alpha'(\omega)}(P)=1$ iff $\wt\omega(P)=1,$
i.e. $\wt\omega(\alpha''(P))=1$ iff $\wt\omega(P)=1.$}
By uniqueness of the projection associated with $\S_0$
(cf. 3.6.11~[Pe]) we conclude that $\alpha''(P)=P.$ Thus
$$\alpha''(\pi\s\l.(\f)'')=\alpha''(P\pi\s\f.(\f)'')=
\pi\s\l.(\f)''=\l''$$
and so $\alpha''\restriction\l''\in\aut\l''.$

For the next step, one would like to prove that $\alpha''$ preserves $\l\subset\l'',$
but here some problems occur. 
Given that $\alpha'\in A(\S(\f))$ preserves $\S_0,$
and that $\S_0$ is the state space of $\l$ (via $\theta$) one would like
to restrict $\alpha'$ to $\S_0,$ and then use Kadison's result that 
there is a bijection between affine w*--continuous homeomorphisms of
$\S(\l)$ and the Jordan automorphism on $\l.$ However, this can fail because
whilst we know that $\alpha'$ is continuous w.r.t. the w*--topology on $\S_0$
coming from $\f,$ when we consider $\S_0$ as the state space of $\l$
we change the w*--topology, and it is not obvious anymore that
$\alpha'\restriction\S_0$ is a homeomorphism w.r.t. this new w*--topology.
We need to know something more of the construction of $\l.$


\thrm Theorem 3.2."
Let $\l$ be an ideal host  for a pair $(\cl F.,\,\S_0)$
and let $\alpha:G\to\aut\f$ be a homomorphism where $G$ is locally compact, and such
that $\alpha''$ preserves $\l.$ If the maps $g\to\omega(\alpha''_g(A))$
are continuous for all $\omega\in\S_0$ and $A\in\l,$ and if $\l$ is separable,
then we have the norm--continuity
${\|\alpha_g''(A)-A\|}\to 0$ as $g\to e$ for all $A\in\l.$"
Apply Theorem 7 of Aarnes [Aa].

This theorem is not really optimal, since one would like to start with conditions
on the original action $\alpha:G\to\aut\f$ only.

\beginsection Acknowledgements.

I am deeply grateful to Detlev Buchholz, who has been as much responsible for the
development of host algebras as I have. Not only did the current line of thought develop out 
of his question at a seminar I gave in 1997
and many subsequent lively discussions, but he also read
several previous misguided attempts, discovering
deeply buried and serious errors.
Moreover, over the years he prodded me to return to
this project, and most recently at G\"ottingen 
provided the opportunity and support for me to pursue it, of which this
paper is the result. Thank you Detlev!

I am also grateful to the Erwin Schr\"odinger Institute in Vienna
(and in particular to Jakob Yngvason and Heide Narnhofer)
who provided me with the opportunity of completing the last part
of the paper.

\beginsection Bibliography.

\item{[Aa]} Aarnes; Continuity of group representations, with applications 
 to C*--algebras, 
J. Funct. Anal. {\bf 5}, 14-36 (1970)
\item{[AP]} C. Akemann, G.K. Pedersen: Complications of
semicontinuity in C*--algebra theory, Duke Math. J.  {\bf 40}, 785--795 (1973)
\item{[BR]} O. Bratteli, D. Robinson,
Operator Algebras and Quantum Statistical Mechanics I. Springer, New York,
1979.
\item{[Di]} J. Dixmier; C*--algebras, North--Holland,Amsterdam, 1977.
\item{[FD]} J.M.G. Fell, R.S. Doran, Representations of *-algebras, 
locally compact groups, and Banach *- algebraic bundles. Vol. 1~\&~2,
Academic Press, 1991.
\item{[GH]} H. Grundling, C.A. Hurst;
Algebraic quantization of systems with a
  gauge degeneracy. Commun. Math. Phys. {\bf 98}, 369--390 (1985)\chop
Grundling, H.: Systems with outer constraints. Gupta--Bleuler
  electromagnetism as an algebraic field theory. Commun. Math. Phys.
 {\bf 114}, 69--91 (1988)
\item{[Gr]} H. Grundling, A Group algebra for inductive limit groups.
Continuity Problems of the Canonical Commutation Relations.
Acta Applicandae Math. {\bf 46}, 107--145, 1997. 
\item{[Ha]}
Haag, R.: Local Quantum Physics. Berlin: Springer Verlag 1992
\item{[HKK]} R. Haag, R.V. Kadison, D. Kastler,
Nets of C*--algebras and Classification of States.
Commun. Math. Phys. {\bf 16}, 81--104 (1970)
\item{[KR]}
Kadison, R.V., Ringrose, J.R.: Fundamentals of the Theory of Operator Algebras II.
  Orlando: Academic Press 1986
\item{[Ku]}
Kusuda, M.: Unique state extension and hereditary C*--subalgebras.
Math. Ann. {\bf 288}, 201--209 (1990)
\item{[Ma]}
Manuceau, J.: C*--alg\`ebre de relations de commutation. Ann.
  Inst. H. Poincar\'e {\bf 8}, 139--161 (1968)
\item{[Mu]}
Murphy, G.J.: C*--Algebras and Operator Theory. Boston:
  Academic Press 1990
\item{[Ni]} Nilsen, M.: C*--bundles and $C_0(X)\hbox{--algebras.}$
Indiana Univ. Math. J., {\bf 45}, 463--477 (1996)
\item{[Pe]} 
Pedersen, G.K.: C*--Algebras and their Automorphism Groups.
  London: Academic Press 1989
\item{[Ri]} M.A. Rieffel: Induced representations of C*--algebras.
Adv. Math. {\bf 13}, 176--257 (1974)
\item{[Se]} I.E. Segal: Representations of the canonical
commutation relations, Carg\'ese lectures in theoretical physics,
Gordon and Breach, 1967.
\item{[Ta]} A. Taylor: Introduction to Functional Analysis.
John Wiley \& Sons, New York 1958.
\item{[Tak]}
Takesaki, M.: Theory of operator algebras I. New York: Springer-Verlag
  1979

\bye